\documentclass[a4paper, 11pt, reqno]{amsart}

\title{Quantum $K$-theoretic Whitney relations for type $C$ flag manifolds}
\author{Takafumi Kouno}
\address{Department of Applied Mathematics, Faculty of Science and Engineering, Waseda University, Japan}
\email{t.kouno@aoni.waseda.jp}

\allowdisplaybreaks[4]

\usepackage[margin=25mm]{geometry}
\usepackage{hyperref}

\usepackage{amssymb}
\usepackage{mathtools}
\mathtoolsset{showonlyrefs}
\usepackage{bm}
\usepackage{dynkin-diagrams}
\usepackage{enumitem}

\numberwithin{equation}{section}

\theoremstyle{plain}
\newtheorem{thm}{Theorem}[section]
\newtheorem{prop}[thm]{Proposition}
\newtheorem{lem}[thm]{Lemma}
\newtheorem{cor}[thm]{Corollary}
\newtheorem*{claim}{Claim}

\newtheorem{ithm}{Theorem}

\theoremstyle{definition}
\newtheorem{defn}[thm]{Definition}

\theoremstyle{remark}
\newtheorem{rem}[thm]{Remark}

\newcommand{\BC}{\mathbb{C}}
\newcommand{\BP}{\mathbb{P}}
\newcommand{\BR}{\mathbb{R}}
\newcommand{\BZ}{\mathbb{Z}}

\newcommand{\udBC}{\underline{\BC}}

\newcommand{\BI}{\mathbf{I}}

\newcommand{\e}{\mathbf{e}}
\newcommand{\bp}{\mathbf{p}}
\newcommand{\bchi}{\bm{\chi}}

\newcommand{\CA}{\mathcal{A}}
\newcommand{\CE}{\mathcal{E}}
\newcommand{\CO}{\mathcal{O}}
\newcommand{\CS}{\mathcal{S}}
\newcommand{\CZ}{\mathcal{Z}}

\newcommand{\SB}{\mathsf{B}}
\newcommand{\SW}{\mathsf{W}}

\newcommand{\st}{\mathsf{t}}

\newcommand{\Fh}{\mathfrak{h}}

\newcommand{\Sp}{\operatorname{Sp}}
\newcommand{\bQG}{\mathbf{Q}_{G}}
\newcommand{\Par}{\overline{\operatorname{Par}}}

\newcommand{\af}{\mathrm{af}}

\newcommand{\ve}{\varepsilon}
\newcommand{\vpi}{\varpi}

\newcommand{\hr}{\widehat{r}}

\newcommand{\vtl}{\vartriangleleft}
\newcommand{\vtr}{\vartriangleright}

\newcommand{\extprod}{\bigwedge\nolimits}
\newcommand{\qprod}{\sideset{}{^{\star}}{\prod}}

\newcommand{\FlC}[1]{\operatorname{Fl}^{C_{#1}}}
\newcommand{\bra}[1]{[\![ #1 ]\!]}
\newcommand{\pra}[1]{(\!( #1 )\!)}
\newcommand{\ti}[1]{\widetilde{#1}}
\newcommand{\pair}[2]{\langle #1, #2 \rangle}

\DeclareMathOperator{\rk}{rk}
\DeclareMathOperator{\Fl}{Fl}
\DeclareMathOperator{\diag}{diag}

\DeclareMathOperator{\gch}{gch}
\DeclareMathOperator{\End}{End}

\DeclareMathOperator{\QBG}{QBG}
\DeclareMathOperator{\init}{init}
\DeclareMathOperator{\ed}{end}
\DeclareMathOperator{\down}{down}
\DeclareMathOperator{\wt}{wt}
\DeclareMathOperator{\sign}{sign}

\newenvironment{enu}{%
\begin{enumerate}[label=\textup{(\arabic*)}, left=0pt]}{%
\end{enumerate}}
\newenvironment{enui}{%
\begin{enumerate}[label=\textup{(\roman*)}, left=0pt]}{%
\end{enumerate}}

\makeatletter
\def\paragraph{\@startsection{paragraph}{4}%
  \z@{0.5\baselineskip}{-\fontdimen2\font}%
  \normalfont}
\def\subparagraph{\@startsection{subparagraph}{5}%
  \z@{0.5\baselineskip}{-\fontdimen2\font}%
  \normalfont}
\makeatother

\begin{document}

\begin{abstract}
We study relations of $\lambda_{y}$-classes associated to tautological bundles over the flag manifold of type $C$ in the quantum $K$-ring. 
These relations are called the quantum $K$-theoretic Whitney relations. 
The strategy of the proof of the quantum $K$-theoretic Whitney relations is based on the method of semi-infinite flag manifolds and the Borel-type presentation. 
In addition, we observe that the quantum $K$-theoretic Whitney relations give a complete set of the defining relations of the quantum $K$-ring. This gives a presentation of the quantum $K$-ring of the flag manifold of type $C$, called the Whitney-type presentation, as a quotient of a polynomial ring, different from the Borel-type presentation. 
\end{abstract}

\maketitle

\section{Introduction}

Let $X$ be a partial flag manifold of arbitrary type. 
Givental \cite{Givental} and Lee \cite{Lee} introduced the \emph{torus-equivariant quantum $K$-ring} $QK_{T}(X)$ of the flag manifold $X$, 
which is a deformation of the ordinary torus-equivariant $K$-ring $K_{T}(X)$; here $T$ denotes the torus. 
The product $\star$ of $QK_{T}(X)$, called the \emph{quantum product}, is defined by using \emph{$K$-theoretic Gromov--Witten invariants}. 
To understand the quantum product, presentations of $QK_{T}(X)$ with generators and relations have been studied, including Borel-type, Coulomb branch, and Whitney presentations. In this paper, we focus on the Whitney presentation. 

For a $T$-equivariant vector bundle $\CE$ over $X$, we denote by $[\CE]$ the class of $\CE$ in $K_{T}(X)$. 
The \emph{Hurzebruch $\lambda_{y}$-class} $\lambda_{y}(\CE) \in K_{T}(X)[y]$ of a $T$-equivariant vector bundle $\CE$, introduced in \cite{Hirzebruch}, is defined as 
\begin{equation}
\lambda_{y}(\CE) := 1 + y [\CE] + y^{2} \left[ \extprod^{2} \CE \right] + \cdots + y^{\rk(\CE)} \left[ \extprod^{\rk(\CE)} \CE \right],  
\end{equation}
where $\rk(\CE)$ denotes the rank of $\CE$. 
For a short exact sequence $0 \rightarrow \CE_{1} \rightarrow \CE_{2} \rightarrow \CE_{3} \rightarrow 0$ of torus-equivariant vector bundles over $X$, we have 
\begin{equation} \label{eq:Whitney}
\lambda_{y}(\CE_{1}) \cdot \lambda_{y}(\CE_{3}) = \lambda_{y}(\CE_{2}) 
\end{equation}
in $K_{T}(X)[y]$, where $\cdot$ denotes the product of $K_{T}(X)$ (or $K_{T}(X)[y]$).
The identity \eqref{eq:Whitney} is called the \emph{$K$-theoretic Whitney relation}. 

Our interest is to quantize the above $K$-theoretic Whitney relations.
More precisely, we study the quantum product $\lambda_{y}(\CE_{1}) \star \lambda_{y}(\CE_{3})$ in $QK_{T}(X)[y]$. 
Let us consider the case $X = \Fl(r_{1}, \ldots, r_{k}; n)$ is a $k$-step flag manifold of type $A$ parametrizing all sequences $(0 \subset V_{1} \subset \cdots \subset V_{k} \subset \BC^{n})$ of vector subspaces of $\BC^{n}$ with $\dim V_{i} = r_{i}$ for $1 \le i \le k$.
Then we can take a tautological filtration $0 = \CS_{0} \subset \CS_{1} \subset \cdots \subset \CS_{k} \subset \udBC^{n}$ of the trivial bundle $\udBC^{n} = X \times \BC^{n}$. 
Gu, Mihalcea, Sharpe, Xu, Zhang, and Zou (\cite{GMSXZZ}) conjectured the result of the quantum product of $\lambda_{y}$-classes associated to tautological vector bundles as follows: 
\begin{equation} \label{eq:QK-Whitney_typeA}
\lambda_{y}(\CS_{j}) \star \lambda_{y}(\CS_{j+1}/\CS_{j}) = \lambda_{y}(\CS_{j+1}) - y^{r_{j+1}-r_{j}} \dfrac{Q_{j}}{1-Q_{j}} \det(\CS_{j+1}/\CS_{j}) \star (\lambda_{y}(\CS_{j}) - \lambda_{y}(\CS_{j-1})) 
\end{equation}
for $1 \le j \le k$, here $Q_{1}, \ldots, Q_{k}$ are Novikov variables of $QK_{T}(X)$, and $\lambda_{y}(\CS_{n+1}) = \lambda_{y}(\udBC^{n}) = \prod_{1 \le j \le n} (1 + yT_{j})$, where $T_{j} \in R(T)$ for $1 \le j \le n$ are given by the decomposition of $\BC^{n}$ into one-dimensional $T$-modules (see \cite[Conjecture~1.1]{GMSXZZ}). 
This identity is called the \emph{quantum $K$-theoretic Whitney relation} (the \emph{$QK$-Whitney relation} for short). 
This $QK$-Whitney relation was proved for Grassmannians $X = \mathrm{Gr}(k, n) := \Fl(k; n)$ by Gu, Mihalcea, Sharpe, and Zou (\cite{GMSZ2}), for $1 \le k \le n$, for the incidence variety $X = \Fl(1, n-1; n)$ by Gu, Mihalcea, Sharpe, Xu, Zhang, and Zou (\cite{GMSXZZ}), and for arbitrary partial flag manifolds in type $A$ by Huq-Kuruvilla (\cite{HK}). 
In addition, Gu, Mihalcea, Sharpe, Xu, Zhang, and Zou showed in \cite[Theorem~1.4]{GMSXZZ} that this conjecture holds under the quantum $K$-theoretic divisor axiom stated by Buch and Mihalcea, the identity for three-pointed $K$-theoretic Gromov--Witten invariants. 
Recently, Lenart, Naito, Sagaki, and Xu proved the quantum $K$-theoretic divisor axiom in \cite{LNSX} (not only for $X = \Fl(r_{1}, \ldots, r_{k}; n)$ but also for arbitrary $X$). 
As a benefit of describing $QK$-Whitney relations, we can obtain a presentation of $QK_{T}(\Fl(r_{1}, \ldots, r_{k}; n))$ as a quotient ring of a polynomial ring. In fact, $QK$-Whitney relations give a complete set of defining relations of $QK_{T}(\Fl(r_{1}, \ldots, r_{k}, n))$ (\cite[Corollary~5.4]{GMSXZZ2}). 

Since the progress of the study of $QK$-Whitney relations is recent, explicit descriptions of $QK$-Whitney relations are known only for type $A$ partial flag manifolds at present; note that conjectures exist beyond type $A$ (for example, \cite{GMSZ}). The purpose of this paper is to describe $QK$-Whitney relations for the flag manifold in type $C$. 
Let $X = \FlC{n}$ be the flag manifold of type $C_{n}$, the manifold parametrizing all sequences 
\begin{equation}
(0 \subset V_{1} \subset \cdots \subset V_{n} \subset \BC^{2n})
\end{equation}
of isotropic subspaces $V_{1}, \ldots, V_{n}$ with respect to a fixed symplectic form with $\dim V_{k} = k$ for $1 \le k \le n$. 
Then we can take a tautological filtration
\begin{equation}
0 = \CS_{0} \subset \CS_{1} \subset \cdots \subset \CS_{n} \subset \udBC^{2n}
\end{equation}
of subbundles of the trivial bundle $\udBC^{2n} = \FlC{n} \times \BC^{2n}$ over $\FlC{n}$ such that $\rk(\CS_{k}) = k$ for $1 \le k \le n$. 
For a $T$-equivariant vector bundle $\CE$ over $\FlC{n}$, we denote by $\CE^{\vee}$ its dual bundle. 
We denote by $Q_{1}, \ldots, Q_{n}$ the Novikov variables of $QK_{T}(\FlC{n})$. We understand that $Q_{0} = 0$ and $\CS_{-1} = 0$. 
Our main results are as follows. 

\begin{ithm}[$=$\,Theorem~\ref{thm:QK-Whitney_A}\,$+$\,Theorem~\ref{thm:QK-Whitney_B}] \label{ithm:QK-Whitney}
\begin{enu}
\item In $QK_{T}(\FlC{n})[y]$, for $1 \le k \le n$, we have 
\begin{equation} \label{eq:QK-Whitney_1_intro}
\lambda_{y}(\CS_{k-1}) \star \lambda_{y}(\CS_{k}/\CS_{k-1}) = \lambda_{y}(\CS_{k}) - y \frac{Q_{k-1}}{1-Q_{k-1}} [\CS_{k}/\CS_{k-1}] \star (\lambda_{y}(\CS_{k-1}) - \lambda_{y}(\CS_{k-2})). 
\end{equation}

\item In $QK_{T}(\FlC{n})[y]$, for $1 \le k \le n$, we have
\begin{equation} \label{eq:QK-Whitney_2_intro}
\lambda_{y}(\CS_{k-1}^{\vee}) \star \lambda_{y}((\CS_{k}/\CS_{k-1})^{\vee}) = \lambda_{y}(\CS_{k}^{\vee}) - y \frac{Q_{k-1}}{1-Q_{k-1}} [(\CS_{k}/\CS_{k-1})^{\vee}] \star (\lambda_{y}(\CS_{k-1}^{\vee}) - \lambda_{y}(\CS_{k-2}^{\vee})). 
\end{equation}

\item In $QK_{T}(\FlC{n})[y]$, we have 
\begin{equation} \label{eq:QK-Whitney_3_intro}
\begin{split}
& \lambda_{y}(\CS_{n}) \star \lambda_{y}(\CS_{n}^{\vee}) = \lambda_{y}(\udBC^{2n}) \\ 
& \qquad - y^{2} \sum_{p=1}^{n} \frac{Q_{p} \cdots Q_{n}}{1-Q_{p-1}} (\lambda_{y}(\CS_{p-1}) - Q_{p-1} \lambda_{y}(\CS_{p-2})) \star (\lambda_{y}(\CS_{p-1}^{\vee}) - Q_{p-1} \lambda_{y}(\CS_{p-2}^{\vee})). 
\end{split}
\end{equation}

\end{enu}
\end{ithm}
While \eqref{eq:QK-Whitney_1_intro} and \eqref{eq:QK-Whitney_2_intro} are analogous to \eqref{eq:QK-Whitney_typeA}, \eqref{eq:QK-Whitney_3_intro} is specific to the type $C$ case. 

Then we observe that relations \eqref{eq:QK-Whitney_1_intro}, \eqref{eq:QK-Whitney_2_intro}, and \eqref{eq:QK-Whitney_3_intro} give a complete set of the defining relations of $QK_{T}(\FlC{n})$ and hence we obtain a presentation of $QK_{T}(\FlC{n})$, which we call the \emph{Whitney-type presentation}. 

\begin{ithm}[$=$\,Theorem~\ref{thm:Whitney_presentation}] \label{ithm:Whitney_presentation}
Define a ring $R_{Q}^{\SW}$ by 
\begin{equation}
R_{Q}^{\SW} := (R(T)\bra{Q_{1}, \ldots, Q_{n}}) [F_{d}^{k}, G_{d}^{k}, x_{j}, y_{j} \mid 1 \le k \le n, \ 1 \le d \le k, \ 1 \le j \le n]. 
\end{equation}
Take an ideal $I_{Q}^{\SW}$ of $R_{Q}^{\SW}$ generated by 
\begin{align}
& F_{d}^{k} - \left( F_{d}^{k-1} + \frac{1}{1-Q_{k-1}} (F_{d-1}^{k-1} - Q_{k-1} F_{d-1}^{k-2})x_{k} \right) \quad \text{for $1 \le k \le n$ and $1 \le d \le k$}, \\ 
& G_{d}^{k} - \left( G_{d}^{k-1} + \frac{1}{1-Q_{k-1}} (G_{d-1}^{k-1} - Q_{k-1} G_{d-1}^{k-2})y_{k} \right) \quad \text{for $1 \le k \le n$ and $1 \le d \le k$}, \\ 
\begin{split}
&\sum_{k = 0}^{d} F_{k}^{n} G_{d-k}^{n} - \Biggl( e_{d}(T_{1}, T_{2}, \ldots, T_{n}, T_{n}^{-1}, \ldots, T_{2}^{-1}, T_{1}^{-1}) \\ 
& \quad - \sum_{p = 1}^{n} \frac{Q_{p} \cdots Q_{n}}{1-Q_{p-1}} \sum_{r=0}^{d-2} (F_{r}^{p-1} - Q_{p-1} F_{r}^{p-2})(G_{d-2-r}^{p-1} - Q_{p-1}G_{d-2-r}^{p-2}) \Biggr) \quad \text{for $1 \le d \le n$}, 
\end{split} \\ 
& F_{1}^{1} - x_{1}, \\
& G_{1}^{1} - y_{1}, \\ 
& x_{j}y_{j} = (1-Q_{j-1})(1-Q_{j}) \quad \text{for $1 \le j \le n$}; 
\end{align}
here, $e_{d}$ denotes the $d$-th elementary symmetric polynomial, $T_{k} \in R(T)$ for $1 \le k \le 2n$ are obtained from the decomposition of $\BC^{2n}$ into one-dimensional $T$-modules; note that $T_{n+1} = T_{n}^{-1}$, $T_{n+2} = T_{n-1}^{-1}, \ldots$, $T_{2n} = T_{1}^{-1}$. 
In addition, we understand that 
$F_{d}^{-1} = F_{d}^{0} = 0$ for $d > 0$, $F_{d}^{k} = 0$ for $1 \le k \le n$ and $d > k$, $F_{0}^{k} = 1$ for $-1 \le k \le n$. 

Then there exists an $R(T)$-algebra isomorphism $R_{Q}^{\SW}/I_{Q}^{\SW} \xrightarrow{\sim} QK_{T}(\FlC{n})$ such that 
\begin{align}
F_{d}^{k} & \mapsto \left[ \extprod^{d} \CS_{k} \right] & & \text{for $1 \le k \le n$ and $1 \le d \le k$}, \\ 
G_{d}^{k} & \mapsto \left[ \extprod^{d} \CS_{k}^{\vee} \right] & & \text{for $1 \le k \le n$ and $1 \le d \le k$}, \\ 
x_{j} & \mapsto [\CS_{j}/\CS_{j-1}] & & \text{for $1 \le j \le n$}, \\ 
y_{j} & \mapsto [(\CS_{j}/\CS_{j-1})^{\vee}] & & \text{for $1 \le j \le n$}. 
\end{align}
\end{ithm}

To prove Theorem~\ref{ithm:QK-Whitney}, we use the theory of semi-infinite flag manifolds. 
Kato (\cite{Kato}) proved that the $K$-group of the semi-infinite flag manifold is isomorphic to the quantum $K$-ring of the flag manifold as $R(T)$-modules. 
Through this isomorphism, we can translate certain quantum products in $QK_{T}(\FlC{n})$ that appear in the computation of the left-hand side of $QK$-Whitney relations into tensor products in the $K$-group of the semi-infinite flag manifold. 
Since tensor products in the $K$-group of the semi-infinite flag manifold satisfy a multiplicative property, we can obtain identities of line bundle classes in $QK_{T}(\FlC{n})$ (Corollary~\ref{cor:multiple_line_bundles}). 
These identities enable us to compute $\lambda_{y}(\CS_{k-1}) \star \lambda_{y}(\CS_{k}/\CS_{k-1})$ and $\lambda_{y}(\CS_{k-1}^{\vee}) \star \lambda_{y}((\CS_{k}/\CS_{k-1})^{\vee})$. 
To consider the product $\lambda_{y}(\CS_{n}) \star \lambda_{y}(\CS_{n}^{\vee})$, we use an additional method, the Borel-type presentation given in \cite{KN2}. This gives a relation in $QK_{T}(\FlC{n})$ (Proposition~\ref{prop:QK-Borel_rel}) which enables us to simplify such a product. 

To prove Theorem~\ref{ithm:Whitney_presentation}, we use Nakayama-type lemma established by Gu, Mihalcea, Sharpe, Xu, Zhang, and Zou (\cite{GMSXZZ2}). This lemma shows that the isomorphism for the classical $K$-ring $K_{T}(X)$ lifts to that for the quantum $K$-ring $QK_{T}(X)$ under conditions (Theorem~\ref{thm:Nakayama}). By applying this lemma, our Whitney-type presentation reduces to the well-known presentation of $K_{T}(\FlC{n})$, and hence we obtain the theorem. 

This paper is organized as follows. In Section~\ref{sec:preliminaries}, we fix notation for root systems; in particular, for the type $C$ root system. 
Also, we briefly review the flag manifold and its $K$-ring. 
In Section~\ref{sec:QK-Whitney}, we give explicit descriptions of $QK$-Whitney relations in $QK_{T}(\FlC{n})$ (Theorem~\ref{ithm:QK-Whitney}) and the Whitney-type presentation of $QK_{T}(\FlC{n})$ (Theorem~\ref{ithm:Whitney_presentation}), which are the main results of this paper. 
In Section~\ref{sec:method}, we review the semi-infinite flag manifold and its $K$-group. Then we compute some quantum products necessary to study $QK$-Whitney relations, where some proofs are deferred to Appendix~\ref{sec:appendix}. 
In Section~\ref{sec:QK-Whitney_computation}, we give a proof of Theorem~\ref{ithm:QK-Whitney}. In the computation of $\lambda_{y}(\CS_{n}) \star \lambda_{y}(\CS_{n}^{\vee})$, we review Borel-type relations of $QK_{T}(\FlC{n})$. 
In Section~\ref{sec:Nakayama}, we review the Nakayama-type lemma and give a proof of Theorem~\ref{ithm:Whitney_presentation}. 

\subsection*{Acknowledgments}
The author is grateful to Takeshi Ikeda and Satoshi Naito for helpful discussions. 
This work was supported by JSPS Grant-in-Aid for Research Activity Start-up 24K22842. 

\section{Preliminaries} \label{sec:preliminaries}

\subsection{Notation for root systems}
First, we fix the notation for root systems. 
Let $G$ be a connected, simply-connected, simple algebraic group over $\BC$. 
Take a maximal torus $T$ of $G$ and a Borel subgroup $B$ of $G$ containing $T$. 
Also, we take a Borel subgroup $B^{-}$ of $G$ opposite to $G$; that is, $B \cap B^{-} = T$.
Let $\Fh$ be the Lie algebra of $T$, and denote by $\pair{\cdot}{\cdot}$ the canonical pairing $\Fh^{\ast} \times \Fh \rightarrow \BC$; here $\Fh^{\ast} = \operatorname{Hom}_{\BC}(\Fh, \BC)$. 

Take the root system $\Delta$ of $G$, and take the set $\Delta^{+}$ of all positive roots. 
We denote by $I$ the index set of the Dynkin diagram of $G$, and by $\alpha_{i}$ for $i \in I$ the simple root corresponding to $i$. 
For $\alpha \in \Delta$, we denote by $\alpha^{\vee}$ the coroot corresponding to $\alpha$. 
Then we define the root lattice $Q$ and the coroot lattice $Q^{\vee}$ of $G$ by 
\begin{equation}
Q := \bigoplus_{i \in I} \BZ \alpha_{i}, \quad Q^{\vee} := \bigoplus_{i \in I} \BZ \alpha_{i}^{\vee}. 
\end{equation}
Also, for $\alpha \in \Delta$, we set 
\begin{equation}
|\alpha| := \begin{cases}
\alpha & \text{if $\alpha \in \Delta^{+}$}, \\ 
-\alpha & \text{if $\alpha \in -\Delta^{+}$}. 
\end{cases}
\end{equation}

Let $W$ be the Weyl group of $G$. Then for $w \in W$, we have a lift $\dot{w} \in G$ of $w$. 
For $i \in I$, we denote by $s_{i}$ the simple reflection corresponding to $i$. 
In addition, we denote by $w_{\circ} \in W$ the longest element of $W$, and by $e$ the identity element of $W$. 

Let $P = \bigoplus_{i \in I} \BZ \vpi_{i}$ be the weight lattice of $G$, where $\vpi_{i}$ for $i \in I$ denotes the $i$-th fundamental weight. 

\subsection{Flag manifolds}

The manifold $G/B$ is called the \emph{flag manifold}.  
We define the (opposite) \emph{Schubert cell} $X^{w, \circ}$ as $X^{w, \circ} := B^{-} \dot{w} B/B$. 
Then we have a decomposition of $G/B$ as follows: 
\begin{equation}
G/B = \bigsqcup_{w \in W} X^{w, \circ}. 
\end{equation}
The Zariski closure $X^{w} := \overline{X^{w, \circ}}$ of the Schubert cell $X^{w, \circ}$ for $w \in W$ is called the (opposite) \emph{Schubert variety}. 

For $\lambda \in P$, we denote by $\BC_{\lambda}$ the one-dimensional representation of $B$ whose weight is $\lambda$. Then for $\lambda \in P$, we define the line bundle $\CO_{G/B}(\lambda)$ over $G/B$ as $G \times_{B} \BC_{-\lambda} \twoheadrightarrow G/B$.

\subsection{The classical and the quantum \texorpdfstring{$K$}{K}-ring of the flag manifold}
We recall the notion of the classical and the quantum $K$-ring of the flag manifold. 
We denote by $R(T)$ the representation ring of $T$. 
Let $K_{T}(G/B)$ be the \emph{$T$-equivariant $K$-group} of $G/B$, which is the Grothendieck group of the category of $T$-equivariant vector bundles over $G/B$. 
It is defined as the $R(T)$-module generated by (formal elements) $[\CE]$ for $T$-equivariant vector bundles $\CE$ over $G/B$, 
with relations $[\CE_{2}] = [\CE_{1}] + [\CE_{3}]$ for all short exact sequences $0 \rightarrow \CE_{1} \rightarrow \CE_{2} \rightarrow \CE_{3} \rightarrow 0$. 
Then $K_{T}(G/B)$ has a ring structure with product $[\CE_{1}] \cdot [\CE_{2}] = [\CE_{1} \otimes \CE_{2}]$ for $T$-equivariant vector bundles $\CE_{1}$, $\CE_{2}$ over $G/B$. 
Since $G/B$ is smooth, $K_{T}(G/B)$ is identical to the Grothendieck ring of the category of $T$-equivariant coherent sheaves over $G/B$. 
For a $T$-equivariant coherent sheaf $\CO$, we denote by $[\CO]$ the class of $\CO$ in $K_{T}(G/B)$. 
For $w \in W$, the class $[\CO^{w}]$ of the Schubert variety $X^{w}$ for $w \in W$ is called a \emph{Schubert class}. 
It is well-known that $\{[\CO^{w}] \mid w \in W\}$ is an $R(T)$-basis of $K_{T}(G/B)$. 

Let $n$ be the rank of $G$. Then the \emph{$T$-equivariant quantum $K$-ring} $QK_{T}(G/B)$ of $G/B$ is defined, as an $R(T)$-module, by
\begin{equation}
QK_{T}(G/B) := K_{T}(G/B) \otimes_{R(T)} R(T) \bra{Q_{1}, \ldots, Q_{n}}. 
\end{equation}
Its product, denoted by $\star$, is defined by using $K$-theoretic Gromov--Witten invariants. For details, see \cite{Givental} and \cite{Lee}. 
It follows that the Schubert classes form an $R(T)\bra{Q_{1}, \ldots, Q_{n}}$-basis of $QK_{T}(G/B)$. 

\subsection{Type \texorpdfstring{$C$}{C} case}
In this section, we concentrate on the case where $G$ is of type $C_{n}$. Note that these settings are the same as those in \cite[\S 2.3]{KN2}. 
Let $J_{2n}^{-}$ be the $2n \times 2n$-matrix 
\begin{equation}
J_{2n}^{-} = \begin{bmatrix} O & J_{n} \\ -J_{n} & O \end{bmatrix}, 
\end{equation}
where $J_{n}$ is the $n \times n$-matrix defined as 
\begin{equation}
J_{n} = \begin{bmatrix}
0 & 0 & \cdots & 0 & 1 \\ 
0 & 0 & \cdots & 1 & 0 \\ 
\vdots & \vdots & & \vdots & \vdots \\ 
0 & 1 & \cdots & 0 & 0 \\ 
1 & 0 & \cdots & 0 & 0 
\end{bmatrix}. 
\end{equation}
We define the symplectic form $(\cdot, \cdot)$ on $\BC^{2n}$ by $(v, w) := {}^{t}v J_{2n}^{-} w$ for $v, w \in \BC^{2n}$. 
We denote by $M_{2n}(\BC)$ the set of all $2n \times 2n$-matrices with complex entries. 
Then we set 
\begin{equation}
\Sp_{2n}(\BC) := \{ A \in M_{2n}(\BC) \mid {}^{t}A J_{2n}^{-} A = J_{2n}^{-} \}. 
\end{equation}
This $\Sp_{2n}(\BC)$ is called the \emph{symplectic group} with respect to $(\cdot, \cdot)$, which is the algebraic group of type $C_{n}$. 
Then we can take a maximal torus $T$ of $G$ as 
\begin{equation}
T = \{ \diag(z_{1}, z_{2}, \ldots, z_{n}, z_{n}^{-1}, \ldots, z_{2}^{-1}, z_{1}^{-1}) \mid z_{1}, \ldots, z_{n} \in \BC^{\ast} = \BC \setminus \{0\} \}, 
\end{equation}
where $\diag(a_{1}, \ldots, a_{2n})$ denotes the $2n \times 2n$-diagonal matrix with entries $a_{1}, \ldots, a_{2n}$. 
The Lie algebra $\Fh$ of $T$ is of the form 
\begin{equation}
\Fh = \{ \diag(a_{1}, a_{2}, \ldots, a_{n}, -a_{n}, \ldots, -a_{2}, -a_{1}) \mid a_{1}, \ldots, a_{n} \in \BC \}. 
\end{equation}
The Dynkin diagram of $\Sp_{2n}(\BC)$ is as follows: 
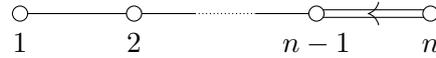
\begin{equation}
\dynkin[labels={1,2,n-1,n}, root radius=0.1cm, edge length=1.5cm, text style/.style={scale=1}]{C}{oo.oo}
\end{equation}
We define $\ve_{j} \in \Fh^{\ast}$ by 
\begin{equation}
\ve_{j}(\diag(a_{1}, a_{2}, \ldots, a_{n}, -a_{n}, \ldots, -a_{2}, -a_{1})) := a_{j}
\end{equation}
for $\diag(a_{1}, a_{2}, \ldots, a_{n}, -a_{n}, \ldots, -a_{2}, -a_{1}) \in \Fh$. 
Then the $j$-th fundamental weight $\vpi_{j}$ for $j \in I := \{1, \ldots, n\}$ is given as 
\begin{equation}
\vpi_{j} = \ve_{1} + \cdots + \ve_{j}
\end{equation}
and the weight lattice $P$ of $\Sp_{2n}(\BC)$ is given as 
\begin{equation}
P = \bigoplus_{1 \le j \le n} \BZ \ve_{j}. 
\end{equation}
Hence we have 
\begin{equation} \label{eq:rep_ring}
R(T) = \BZ[P] = \BZ[\e^{\pm \ve_{1}}, \ldots, \e^{\pm \ve_{n}}]. 
\end{equation}

The action of the Weyl group $W$ on $P$ can be regarded as a signed-permutation on the set $[1, \overline{1}] := \{1 < 2 < \cdots < n < \overline{n} < \cdots < \overline{2} < \overline{1}\}$ via $w(\ve_{j}) = \ve_{w(j)}$ for $w \in W$ and $1 \le j \le \overline{1}$; here we set $\ve_{\overline{j}} := -\ve_{j}$ for $1 \le j \le n$. 
In this paper, we also use the notation of a subinterval $[a, b]$ of $[1, \overline{1}]$ for $1 \le a < b \le \overline{1}$, defined as 
\begin{equation}
[a, b] := \{ j \in [1, \overline{1}] \mid a \le j \le b \}. 
\end{equation}

Let $B$ be the Borel subgroup of $\Sp_{2n}(\BC)$ composed of all upper triangular matrices of $\Sp_{2n}(\BC)$. 
Then the flag manifold $\FlC{n} := \Sp_{2n}(\BC)/B$ of type $C_{n}$ is a manifold parametrizing all sequences
\begin{equation}
(0 \subset V_{1} \subset V_{2} \subset \cdots \subset V_{n} \subset \BC^{2n})
\end{equation}
of isotropic subspaces of $\BC^{2n}$ with respect to the symplectic form $(\cdot, \cdot)$ such that $\dim V_{k} = k$ for $1 \le k \le n$. 
In this situation, we can take the tautological filtration
\begin{equation}
0 = \CS_{0} \subset \CS_{1} \subset \CS_{2} \subset \cdots \subset \CS_{n} \subset \udBC^{2n}
\end{equation}
by subbundles of the trivial bundle $\udBC^{2n} = \FlC{n} \times \BC^{2n}$ over $\FlC{n}$ such that $\rk(\CS_{k}) = k$ for $0 \le k \le n$. 
Then we have the sequence
\begin{equation}
\udBC^{2n} \twoheadrightarrow \CS_{n}^{\vee} \twoheadrightarrow \cdots \twoheadrightarrow \CS_{2}^{\vee} \twoheadrightarrow \CS_{1}^{\vee} \twoheadrightarrow \CS_{0}^{\vee} = 0
\end{equation}
of dual bundles. 

\begin{rem}
For $J \subset \{1, \ldots, n\}$, set 
\begin{equation}
\ve_{J} := \sum_{j \in J} \ve_{j} = \sum_{j \in L_{J}} \vpi_{j} - \sum_{j \in M_{J}} \vpi_{j};  
\end{equation}
where 
\begin{equation}
L_{J} := \{j \in [1, n] \mid \text{$j \in J$ and $j+1 \notin J$}\}, \quad M_{J} := \{j \in [1,n] \mid \text{$j \notin J$ and $j+1 \in J$}\}. 
\end{equation}
In $K_{T}(\FlC{n})$, for $1 \le k \le n$, it follows that 
\begin{align}
\left[ \extprod^{d} \CS_{k} \right] &= \sum_{\substack{J \subset [1, k] \\ |J| = d}} [\CO_{G/B} (-\ve_{J})], \\ 
\left[ \extprod^{d} \CS_{k}^{\vee} \right] &= \sum_{\substack{J \subset [1, k] \\ |J| = d}} [\CO_{G/B} (\ve_{J})], \\ 
[\CS_{k}/\CS_{k-1}] &= [\CO_{G/B}(-\ve_{k})], \\ 
[(\CS_{k}/\CS_{k-1})^{\vee}] &= [\CO_{G/B}(\ve_{k})]. 
\end{align}
In addition, in $K_{T}(\FlC{n})$, for $0 \le d \le 2n$, we see that
\begin{equation}
\left[ \extprod^{d} \udBC^{2n} \right] = e_{d}(\e^{\ve_{1}}, \e^{\ve_{2}}, \ldots, \e^{\ve_{n}}, \e^{-\ve_{n}}, \ldots, \e^{-\ve_{2}}, \e^{-\ve_{1}}), 
\end{equation}
where $e_{d}(x_{1}, \ldots, x_{m})$ denotes the $d$-th elementary symmetric polynomial in variables $x_{1}, \ldots, x_{m}$. 
\end{rem}

\subsection{\texorpdfstring{$K$}{K}-theoretic Whitney relations}
Let $y$ be a formal variable. For a $T$-equivariant vector bundle $\CE$ over $G/B$ with rank $r$, the \emph{Hirzebruch $\lambda_{y}$-class} $\lambda_{y}(\CE)$, introduced in \cite{Hirzebruch}, is defined as 
\begin{equation}
\lambda_{y}(\CE) := 1 + y [\CE] + y^{2} \left[ \extprod^{2} \CE \right] + \cdots + y^{r} \left[ \extprod^{r} \CE \right] \in K_{T}(G/B)[y]. 
\end{equation}
For a short exact sequence $0 \rightarrow \CE_{1} \rightarrow \CE_{2} \rightarrow \CE_{3} \rightarrow 0$ of $T$-equivariant vector bundles over $G/B$, it is known that $\lambda_{y}(\CE_{2}) = \lambda_{y}(\CE_{1}) \cdot \lambda_{y}(\CE_{3})$. 
We call this identity a \emph{$K$-theoretic Whitney relation} (\emph{$K$-Whitney relation} for short). 
In the case $G = \Sp_{2n}(\BC)$, 
short exact sequences
\begin{gather}
0 \rightarrow \CS_{k-1} \rightarrow \CS_{k} \rightarrow \CS_{k}/\CS_{k-1} \rightarrow 0, \\ 
0 \rightarrow (\CS_{k}/\CS_{k-1})^{\vee} \rightarrow \CS_{k}^{\vee} \rightarrow \CS_{k-1}^{\vee} \rightarrow 0, \\ 
0 \rightarrow \CS_{n} \rightarrow \udBC^{2n} \rightarrow \CS_{n}^{\vee} \rightarrow 0 
\end{gather}
for $1 \le k \le n$, imply
\begin{align}
\lambda_{y}(\CS_{k}) &= \lambda_{y}(\CS_{k-1}) \cdot \lambda_{y}(\CS_{k}/\CS_{k-1}), \label{eq:Whitney_classical_1} \\ 
\lambda_{y}(\CS_{k}^{\vee}) &= \lambda_{y}(\CS_{k-1}^{\vee}) \cdot \lambda_{y}((\CS_{k}/\CS_{k-1})^{\vee}), \label{eq:Whitney_classical_2} \\ 
\lambda_{y}(\udBC^{2n}) &= \lambda_{y}(\CS_{n}) \cdot \lambda_{y}(\CS_{n}^{\vee}), \label{eq:Whitney_classical_3}
\end{align}
respectively. In the next section, we consider the relations analogous to \eqref{eq:Whitney_classical_1}, \eqref{eq:Whitney_classical_2}, and \eqref{eq:Whitney_classical_3} in the quantum $K$-theory, which are the main results of this paper.

\section{Whitney-type presentation of the quantum \texorpdfstring{$K$}{K}-ring} \label{sec:QK-Whitney}
We describe the relations in $QK_{T}(\FlC{n})[y]$ analogous to \eqref{eq:Whitney_classical_1}, \eqref{eq:Whitney_classical_2}, and \eqref{eq:Whitney_classical_3}. 
For this purpose, we enumerate the results of the products 
\begin{align}
& \lambda_{y}(\CS_{k-1}) \star \lambda_{y}(\CS_{k}/\CS_{k-1}) & & \text{for $1 \le k \le n$}, \\ 
& \lambda_{y}(\CS_{k-1}^{\vee}) \star \lambda_{y}((\CS_{k}/\CS_{k-1})^{\vee}) & & \text{for $1 \le k \le n$, and} \\ 
&\lambda_{y}(\CS_{n}) \star \lambda_{y}(\CS_{n}^{\vee}). & & 
\end{align} 
We can deduce that these relations form a system of defining relations for $QK_{T}(\FlC{n})$, and hence we obtain a presentation of $QK_{T}(\FlC{n})$. 

\subsection{Quantum \texorpdfstring{$K$}{K}-theoretic Whitney relations}
Throughout this paper, we understand that $Q_{0} = 0$ and $\CS_{-1} = 0$. 
First, the following two identities are relations analogous to \eqref{eq:Whitney_classical_1} and \eqref{eq:Whitney_classical_2}. 
\begin{thm} \label{thm:QK-Whitney_A}
\begin{enu}
\item In $QK_{T}(\FlC{n})[y]$, for $1 \le k \le n$, we have 
\begin{equation} \label{eq:QK-Whitney_1}
\lambda_{y}(\CS_{k-1}) \star \lambda_{y}(\CS_{k}/\CS_{k-1}) = \lambda_{y}(\CS_{k}) - y \frac{Q_{k-1}}{1-Q_{k-1}} [\CS_{k}/\CS_{k-1}] \star (\lambda_{y}(\CS_{k-1}) - \lambda_{y}(\CS_{k-2})). 
\end{equation}

\item In $QK_{T}(\FlC{n})[y]$, for $1 \le k \le n$, we have
\begin{equation} \label{eq:QK-Whitney_2}
\lambda_{y}(\CS_{k-1}^{\vee}) \star \lambda_{y}((\CS_{k}/\CS_{k-1})^{\vee}) = \lambda_{y}(\CS_{k}^{\vee}) - y \frac{Q_{k-1}}{1-Q_{k-1}} [(\CS_{k}/\CS_{k-1})^{\vee}] \star (\lambda_{y}(\CS_{k-1}^{\vee}) - \lambda_{y}(\CS_{k-2}^{\vee})). 
\end{equation}
\end{enu}
\end{thm}

Note that the above identities are similar to those for the type $A$ flag manifold; see \eqref{eq:QK-Whitney_typeA}. 
On the other hand, the following identity, analogous to \eqref{eq:Whitney_classical_3}, is specific to the type $C$ case. 
\begin{thm} \label{thm:QK-Whitney_B}
In $QK_{T}(\FlC{n})[y]$, we have 
\begin{equation} \label{eq:QK-Whitney_3}
\begin{split}
& \lambda_{y}(\CS_{n}) \star \lambda_{y}(\CS_{n}^{\vee}) = \lambda_{y}(\udBC^{2n}) \\ 
& \qquad - y^{2} \sum_{p=1}^{n} \frac{Q_{p} \cdots Q_{n}}{1-Q_{p-1}} (\lambda_{y}(\CS_{p-1}) - Q_{p-1} \lambda_{y}(\CS_{p-2})) \star (\lambda_{y}(\CS_{p-1}^{\vee}) - Q_{p-1} \lambda_{y}(\CS_{p-2}^{\vee})). 
\end{split}
\end{equation}
\end{thm}
We call identities \eqref{eq:QK-Whitney_1}, \eqref{eq:QK-Whitney_2}, and \eqref{eq:QK-Whitney_3} \emph{quantum $K$-theoretic Whitney relations} (\emph{$QK$-Whitney relations} for short) for $\FlC{n}$. 

\subsection{The Whitney-type presentation of the quantum \texorpdfstring{$K$}{K}-ring}
Define a ring $R_{Q}^{\SW}$ by 
\begin{equation}
R_{Q}^{\SW} := (R(T)\bra{Q_{1}, \ldots, Q_{n}}) [F_{d}^{k}, G_{d}^{k}, x_{j}, y_{j} \mid 1 \le k \le n, \ 1 \le d \le k, \ 1 \le j \le n]. 
\end{equation}
Take an ideal $I_{Q}^{\SW}$ of $R_{Q}^{\SW}$ generated by 
\begin{align}
& F_{d}^{k} - \left( F_{d}^{k-1} + \frac{1}{1-Q_{k-1}} (F_{d-1}^{k-1} - Q_{k-1} F_{d-1}^{k-2})x_{k} \right) \quad \text{for $1 \le k \le n$ and $1 \le d \le k$}, \label{eq:QK-Whitney_def1} \\ 
& G_{d}^{k} - \left( G_{d}^{k-1} + \frac{1}{1-Q_{k-1}} (G_{d-1}^{k-1} - Q_{k-1} G_{d-1}^{k-2})y_{k} \right) \quad \text{for $1 \le k \le n$ and $1 \le d \le k$}, \label{eq:QK-Whitney_def2} \\ 
\begin{split}
&\sum_{k = 0}^{d} F_{k}^{n} G_{d-k}^{n} - \Biggl( e_{d}(\e^{\ve_{1}}, \e^{\ve_{2}}, \ldots, \e^{\ve_{n}}, \e^{-\ve_{n}}, \ldots, \e^{-\ve_{2}}, \e^{-\ve_{1}}) \\ 
& \quad - \sum_{p = 1}^{n} \frac{Q_{p} \cdots Q_{n}}{1-Q_{p-1}} \sum_{r=0}^{d-2} (F_{r}^{p-1} - Q_{p-1} F_{r}^{p-2})(G_{d-2-r}^{p-1} - Q_{p-1}G_{d-2-r}^{p-2}) \Biggr) \quad \text{for $1 \le d \le n$} \label{eq:QK-Whitney_def3} , 
\end{split} \\ 
& F_{1}^{1} - x_{1}, \label{eq:initial_1} \\
& G_{1}^{1} - y_{1}, \label{eq:initial_2} \\ 
& x_{j}y_{j} = (1-Q_{j-1})(1-Q_{j}) \quad \text{for $1 \le j \le n$} \label{eq:inverse}; 
\end{align}
here, $e_{d}$ denotes the $d$-th elementary symmetric polynomial and we understand that 
$F_{d}^{-1} = F_{d}^{0} = 0$ for $d \in \BZ_{> 0}$, $F_{d}^{k} = 0$ for $1 \le k \le n$ and $d > k$, $F_{0}^{k} = 1$ for $-1 \le k \le n$. 
Then we can deduce the following $R(T)$-algebra isomorphism, called \emph{Whitney-type presentation} of $QK_{T}(\FlC{n})$. 

\begin{thm} \label{thm:Whitney_presentation}
There exists an $R(T)$-algebra isomorphism $R_{Q}^{\SW}/I_{Q}^{\SW} \xrightarrow{\sim} QK_{T}(\FlC{n})$ such that 
\begin{align}
F_{d}^{k} & \mapsto \left[ \extprod^{d} \CS_{k} \right] & & \text{for $1 \le k \le n$ and $1 \le d \le k$}, \\ 
G_{d}^{k} & \mapsto \left[ \extprod^{d} \CS_{k}^{\vee} \right] & & \text{for $1 \le k \le n$ and $1 \le d \le k$}, \\ 
x_{j} & \mapsto [\CS_{j}/\CS_{j-1}] & & \text{for $1 \le j \le n$}, \\ 
y_{j} & \mapsto [(\CS_{j}/\CS_{j-1})^{\vee}] & & \text{for $1 \le j \le n$}. 
\end{align}
\end{thm}

\section{Computation of quantum products} \label{sec:method}
The purpose of this section is to review an approach to computing quantum products. 
In this computation, the $K$-group of the semi-infinite flag manifold plays a crucial role. 
Thus, in this section, we review facts about the $K$-group of the semi-infinite flag manifold. 
Then we derive formulas for quantum products of line bundles in $QK_{T}(\FlC{n})$ to obtain our desired $QK$-Whitney relations. 
For a while, we assume that $G$ is of an arbitrary type. 

\subsection{The semi-infinite flag manifold and its \texorpdfstring{$K$}{K}-group}
We denote by $N$ the unipotent radical of $B$; note that $B = TN$. Following \cite[\S 3.1]{MNS} (which is based on \cite[\S 1.4 and \S 1.5]{Kato}), \cite[\S 2.1 and \S 2.3]{Orr}, \cite[\S 4.1]{KN2}, and \cite[\S 4.2]{KNS}, we briefly review the definition of the semi-infinite flag manifold $\bQG$. 
For $\lambda \in P$, we denote by $V(\lambda)$ the irreducible highest weight representation of $G$ of weight $\lambda$. 
Then we define $\bQG$ the set of tuples $(\ell_{\lambda})_{\lambda \in P^{+}}$ of lines $\ell_{\lambda} \in \BP(V(\lambda) \otimes_{\BC} \BC\bra{z})$ for $\lambda \in P^{+}$, such that for $\lambda, \mu \in P^{+}$, $\ell_{\lambda + \mu}$ is mapped to $\ell_{\lambda} \otimes_{\BC} \ell_{\mu}$ under the map 
\begin{equation}
V(\lambda + \mu) \otimes_{\BC} \BC\bra{z} \rightarrow (V(\lambda) \otimes_{\BC} \BC\bra{z}) \otimes_{\BC\bra{z}} (V(\mu) \otimes_{\BC} \BC\bra{z})
\end{equation}
induced from the embedding
\begin{equation}
V(\lambda + \mu) \hookrightarrow V(\lambda) \otimes_{\BC} V(\mu)
\end{equation}
of representations of $G$. 
Note that $\ell_{\lambda} \otimes_{\BC} \ell_{\mu}$ above is understood to be its image under the map 
\begin{equation}
(V(\lambda) \otimes_{\BC} \BC\bra{z}) \otimes_{\BC} (V(\mu) \otimes_{\BC} \BC\bra{z}) \rightarrow (V(\lambda) \otimes_{\BC} \BC\bra{z}) \otimes_{\BC\bra{z}} (V(\mu) \otimes_{\BC} \BC\bra{z}). 
\end{equation} 
Then the tuple $(\ell_{\lambda})_{\lambda \in P^{+}}$ is determined by lines $\ell_{\vpi_{i}}$ corresponding to fundamental weights $\vpi_{i}$ for $i \in I$. 
Hence, there exists an embedding
\begin{equation}
\bQG \hookrightarrow \prod_{i \in I} \BP(V(\vpi_{i}) \otimes_{\BC} \BC\bra{z}). \label{eq:closed_embedding}
\end{equation}
Then $\bQG$ is equipped with a scheme structure such that \eqref{eq:closed_embedding} is a closed embedding (see \cite[Lemma~4.4]{KNS}). 
We call $\bQG$ the \emph{semi-infinite flag manifold}. 
Note that $\bQG$ has a $(T \times \BC^{\ast})$-action. 

For $\lambda = \sum_{i \in I} m_{i} \vpi_{i} \in P$ with $m_{i} \in \BZ$, we define the line bundle $\CO_{\bQG}(\lambda)$ over $\bQG$ as the pull-back of $\boxtimes_{i \in I} \CO(m_{i})$ under \eqref{eq:closed_embedding}, where $i$-th $\CO(m_{i})$ is over $\BP(V(\vpi_{i}))$.  
Note that for $\lambda, \mu \in P$, we have $\CO_{\bQG}(\lambda) \otimes \CO_{\bQG}(\mu) = \CO_{\bQG}(\lambda + \mu)$. 

The semi-infinite flag manifold has semi-infinite Schubert subvarieties analogous to Schubert varieties of the flag manifold.
Consider the map $G(\BC\bra{z}) \rightarrow G$ defined as $z \mapsto 0$ and denote by $\BI$ the preimage of $B$ under this map. 
This $\BI$ is called the \emph{Iwahori subgroup} of $G(\BC\bra{z})$. 
Let $W_{\af}$ be the affine Weyl group of $G$. 
For $\xi \in Q^{\vee}$, we denote by $t_{\xi} \in W_{\af}$ the translation element corresponding to $\xi$. 
Then we have $W_{\af} = \{ wt_{\xi} \mid w \in W, \ \xi \in Q^{\vee} \} \simeq W \ltimes Q^{\vee}$. 
Set 
\begin{equation}
Q^{\vee, +} := \bigoplus_{i \in I} \BZ_{\ge 0} \alpha_{i}^{\vee}
\end{equation}
and define $W_{\af}^{\ge 0} := \{wt_{\xi} \mid w \in W, \ \xi \in Q^{\vee, +} \} \simeq W \times Q^{\vee, +}$. 
For $x = wt_{\xi} \in W_{\af}^{\ge 0}$ with $w \in W$ and $\xi \in Q^{\vee, +}$, set $p_{x} := (z^{-\pair{\lambda}{w_{\circ}\xi}} V(\lambda)_{ww_{\circ}\lambda})_{\lambda \in P^{+}}$, where $V(\lambda)_{\mu}$ for $\mu \in P$ is a $\mu$-weight space of $V(\lambda)$. 
Then the set of $(T \times \BC^{\ast})$-fixed points of $\bQG$ is $\{p_{x} \mid x \in W_{\af}^{\ge 0}\}$; that is, $(T \times \BC^{\ast})$-fixed points of $\bQG$ are parametrized by $W_{\af}^{\ge 0}$ (\cite[\S 4.2]{KNS}, \cite[\S 2.3]{Orr}). 
For $x \in W_{\af}^{\ge 0}$, we denote by $\bQG(x)$ the closure of the $\BI$-orbit of $p_{x}$ and we call it the \emph{semi-infinite Schubert variety} associated to $x$. Note that $\bQG(e) = \bQG$. Also, for $x \in W_{\af}^{\ge 0}$, $\bQG(x)$ is a (reduced) closed subscheme of $\bQG$. For $x \in W_{\af}^{\ge 0}$, we denote by $\CO_{\bQG(x)}$ the structure sheaf of $\bQG(x)$. 

Then we define the $K$-group of $\bQG$. 
For $a = \sum_{\lambda \in P} \sum_{k \in \BZ} c_{\lambda, k} q^{k} \e^{\lambda} \in \BZ[q, q^{-1}][P]$ with $c_{\lambda, k} \in \BZ$, we set $|a| := \sum_{\lambda \in P} \sum_{k \in \BZ} |c_{\lambda, k}| q^{k} \e^{\lambda}$. 
Let $\ti{K}'(\bQG)$ be the $\BZ[q, q^{-1}][P]$-module consisting of all formal sums $\sum_{x \in W_{\af}^{\ge 0}} a_{x} [\CO_{\bQG(x)}]$ with $a_{x} \in \BZ[q, q^{-1}][P]$ such that 
\begin{equation}
\sum_{x \in W_{\af}^{\ge 0}} |a_{x}| \gch H^{0}(\bQG, \CO_{\bQG(x)}(\lambda)) \in \BZ[P]\pra{q^{-1}}
\end{equation}
for $\lambda \in P^{++} := \bigoplus_{i \in I} \BZ_{>0} \vpi_{i}$, where $\CO_{\bQG(x)}(\lambda) := \CO_{\bQG(x)} \otimes \CO_{\bQG}(\lambda)$   for $x \in W_{\af}^{\ge 0}$ and $\lambda \in P$, and $\gch$ denotes the character of $(T \times \BC^{\ast})$-weight modules where the character of loop rotation action $\BC^{\ast} \curvearrowright \BC\bra{z}$ is denoted by $q$. 
By \cite[Theorem~1.25]{Kato}, $\ti{K}'(\bQG)$ is equipped with the classes $[\CO_{\bQG(x)}(\lambda)]$ for $x \in W_{\af}^{\ge 0}$ and $\lambda \in P$. 
In particular, $\ti{K}'(\bQG)$ is equipped with the line bundle classes $[\CO_{\bQG}(\lambda)]$ for $\lambda \in P$. 

Let $K_{T \times \BC^{\ast}}(\bQG)$ be the $\BZ[q, q^{-1}][P]$-submodule of $\ti{K}'(\bQG)$ consisting of all formal sums $\sum_{x \in W_{\af}^{\ge 0}} a_{x} [\CO_{\bQG(x)}]$ with $a_{x} \in \BZ[q, q^{-1}][P]$ such that 
\begin{equation}
\sum_{x \in W_{\af}^{\ge 0}} |a_{x}| \in \BZ[P]\pra{q^{-1}}. 
\end{equation}
Then $K_{T \times \BC^{\ast}}(\bQG)$ is in fact a $\BZ[q, q^{-1}][P]$-submodule of $\ti{K}'(\bQG)$ (see \cite[Corollary~4.31]{KNS}). 
Also, for $x \in W_{\af}^{\ge 0}$ and $\lambda \in P$, we have $[\CO_{\bQG(x)}(\lambda)] \in K_{T \times \BC^{\ast}}(\bQG)$ (due to \cite[Theorem~5.16]{KLN} and (the proof of) \cite[Corollary~5.12]{KNS}). 

Let $K_{T}(\bQG)$ be the specialization of $K_{T \times \BC^{\ast}}(\bQG)$ at $q = 1$, and call it the \emph{$T$-equivariant $K$-group} of $\bQG$. 
Then we see that $K_{T}(\bQG)$ is the $R(T) (=\BZ[P])$-module consisting of all infinite linear combinations of $\{ [\CO_{\bQG(x)}] \mid x \in W_{\af}^{\ge 0}\}$ with coefficients in $R(T)$. 
In particular, $\{ [\CO_{\bQG(x)}] \mid x \in W_{\af}^{\ge 0}\}$ forms a topological basis of $K_{T}(\bQG)$ over $R(T)$; that is, we can write each element in $K_{T}(\bQG)$ uniquely as an infinite linear combination of $[\CO_{\bQG(x)}]$, $x \in W_{\af}^{\ge 0}$, with coefficients in $R(T)$ (\cite[Lemma~1.22]{Kato}). 
In addition, we have $[\CO_{\bQG(x)}(\lambda)] \in K_{T}(\bQG)$ for $x \in W_{\af}^{\ge 0}$ and $\lambda \in P$ (\cite[Theorem~1.26]{Kato}); in particular, $[\CO_{\bQG}(\lambda)] = [\CO_{\bQG(e)}(\lambda)] \in K_{T}(\bQG)$ for $\lambda \in P$. 

For $\xi \in Q^{\vee, +}$, we define $\st_{\xi} \in \End(K_{T}(\bQG))$ by $\st_{\xi}([\CO_{\bQG(x)}]) := [\CO_{\bQG(xt_{\xi})}]$ and extend it linearly. If $\xi = \alpha_{i}^{\vee}$ for $i \in I$, then we denote by $\st_{i} := \st_{\alpha_{i}^{\vee}}$. 

We recall a relationship between $QK_{T}(G/B)$ and $K_{T}(\bQG)$. 

\begin{thm}[{\cite[Corollary~3.13 and Theorem~4.17]{Kato}}] \label{thm:QK_vs_semi-infinite}
There exists an isomorphism of $R(T)$-modules
\begin{equation}
\Phi: QK_{T}(G/B) \xrightarrow{\sim} K_{T}(\bQG)
\end{equation}
such that 
\begin{enui}
\item $\Phi(\e^{\mu} Q^{\xi} [\CO^{w}]) = \e^{-\mu} \st_{\xi} [\CO_{\bQG(w)}]$ for $\mu \in P$, $\xi \in Q^{\vee, +}$, and $w \in W$, and 
\item $\Phi([\CO_{G/B}(-\vpi_{i})] \star \CZ) = [\CO_{\bQG}(w_{\circ}\vpi_{i})] \otimes \Phi(\CZ)$ for $i \in I$ and $\CZ \in QK_{T}(G/B)$.
\end{enui} 
\end{thm}

Note that we need $w_{\circ}$ in the right-hand side of the identity in Theorem~\ref{thm:QK_vs_semi-infinite}\,(ii) because the convention of line bundles in this paper is different from that in \cite{Kato} by a twist coming from the involution $-w_{\circ}$.

\subsection{The quantum alcove model}
Lenart and Postnikov (\cite{LP}) introduced the alcove model to describe the Chevalley formula in the classical $K$-theory of the flag manifold. 
The Chevalley formula in $K_{T}(G/B)$ is an expansion identity of the product of a class of line bundle and a Schubert class into a linear combination of Schubert classes; that is, an identity of the form 
\begin{equation}
[\CO_{G/B}(\lambda)] \cdot [\CO^{w}] = \sum_{v \in W} c_{\lambda, w}^{v} [\CO^{v}]
\end{equation}
for $\lambda \in P$ and $w \in W$, where $c_{\lambda, w}^{v} \in R(T)$. 
Then Lenart and Lubovsky (\cite{LL}) introduced the quantum analogue of the alcove model, called the quantum alcove model, to construct a combinatorial model of (tensor product of) specific Kirillov--Reshetikhin crystals. 
Recently, Lenart, Naito, and Sagaki (\cite{LNS}) used this quantum alcove model to describe the Chevalley formula in the $K$-group of the semi-infinite flag manifold. 
In this situation, the Chevalley formula for the semi-infinite case is an expansion identity of a semi-infinite Schubert class twisted by a class of a line bundle into an infinite linear combination of semi-infinite Schubert classes; that is, an identity of the form
\begin{equation}
[\CO_{\bQG(x)}(\lambda)] = \sum_{y \in W_{\af}} d_{\lambda, x}^{y} [\CO_{\bQG(y)}]
\end{equation}
for $\lambda \in P$ and $x \in W_{\af}$, where $d_{\lambda, x}^{y} \in R(T)$. 
In this section, we review the definition of combinatorial objects used in the theory of the quantum alcove model, following \cite[\S 3.2]{LNS}. 
Then we review the combinatorial description of the Chevalley formula in $K_{T}(\bQG)$ and $K_{T}(G/B)$ in the next section. 

To introduce the quantum alcove model, we define the quantum Bruhat graph.
\begin{defn}[{\cite[Definition~6.1]{BFP}}]
The \emph{quantum Bruhat graph} $\QBG(W)$ is a $\Delta^{+}$-labeled directed graph whose vertex set is $W$ and for $x, y \in W$ and $\alpha \in \Delta^{+}$, there exists an edge $x \xrightarrow{\alpha} y$ if one of the following holds: 
\begin{itemize}
\item[(B)] $\ell(y) = \ell(x) + 1$; or 
\item[(Q)] $\ell(y) = \ell(x) - 2 \pair{\rho}{\alpha^{\vee}} + 1$, 
\end{itemize}
where $\ell (\cdot)$ denotes the length function $W \rightarrow \BZ_{\ge 0}$, and $\rho = (1/2) \sum_{\alpha \in \Delta^{+}} \alpha$. 
If (B) (resp., (Q)) holds, then the edge is called a \emph{Bruhat edge} (resp., \emph{quantum edge}). 
\end{defn}

Also, we introduce the notion of alcoves. 
For $\alpha \in \Delta$ and $k \in \BZ$, we define the hyperplane $H_{\alpha, k}$ in $\Fh_{\BR}^{\ast} := P \otimes_{\BZ} \BR$ as 
\begin{equation}
H_{\alpha, k} := \{ \xi \in \Fh_{\BR}^{\ast} \mid \pair{\xi}{\alpha^{\vee}} = k \}. 
\end{equation}
Then an \emph{alcove} is a connected component of 
\begin{equation}
\Fh_{\BR}^{\ast} \setminus \bigcup_{\alpha \in \Delta, \ k \in \BZ} H_{\alpha, k}. 
\end{equation}
In particular, we define the \emph{fundamental alcove} $A_{\circ}$ as 
\begin{equation}
A_{\circ} = \{ \xi \in \Fh_{\BR}^{\ast} \mid \text{$0 < \pair{\xi}{\alpha^{\vee}} < 1$ for all $\alpha \in \Delta^{+}$} \}. 
\end{equation}
For $\lambda \in P$, we set 
\begin{equation}
A_{\lambda} = A_{\circ} + \lambda := \{ \xi + \lambda \mid \xi \in A_{\circ} \}. 
\end{equation}
Then $A_{\lambda}$ is also an alcove. 

We say $A$ and $B$ are \emph{adjacent} if $A$ and $B$ have a common wall. 
For adjacent alcoves $A$ and $B$ and $\alpha \in \Delta$, we write $A \xrightarrow{\alpha} B$ if the common wall of $A$ and $B$ is contained in $H_{\alpha, k}$ for some $k \in \BZ$ and $\alpha$ points the direction from $A$ to $B$. 

\begin{defn}[{\cite[Definitions~5.2 and 5.4]{LP}}]
\begin{enu}
\item A sequence $(A_{0}, \ldots, A_{r})$ of alcoves is called an \emph{alcove path} if $A_{k-1}$ and $A_{k}$ are adjacent for all $1 \le k \le r$. 
We call an alcove path $\Gamma = (A_{0}, \ldots, A_{r})$ \emph{reduced} if the length $r$ of the alcove path $\Gamma$ is minimal among all alcove paths from $A_{0}$ and $A_{r}$. 
\item For $\lambda \in P$, a sequence $(\gamma_{1}, \ldots, \gamma_{r})$ of roots is called a \emph{$\lambda$-chain} if there exists an alcove path $(A_{0}, \ldots, A_{r})$ with $A_{0} = A_{\circ}$ and $A_{r} = A_{-\lambda}$ such that 
\begin{equation}
A_{\circ} = A_{0} \xrightarrow{-\gamma_{1}} A_{1} \xrightarrow{-\gamma_{2}} \cdots \xrightarrow{-\gamma_{r}} A_{r} = A_{-\lambda}. 
\end{equation}
We say that a $\lambda$-chain $(\gamma_{1}, \ldots. \gamma_{r})$ is \emph{reduced} if the corresponding alcove path $(A_{0}, \ldots, A_{r})$ is reduced. 
\end{enu}
\end{defn}

Finally, we introduce admissible subsets, the main objects in the theory of the quantum alcove model. 
\begin{defn}[{\cite[Definition~3.4]{LL}, \cite[Definition~17]{LNS}}]
Let $\lambda \in P$ and $w \in W$. Take a $\lambda$-chain $\Gamma = (\gamma_{1}, \ldots, \gamma_{r})$ (not necessarily reduced). 
Then a set $A = \{i_{1} < \ldots < i_{t}\} \subset \{1, \ldots, r\}$ (possibly empty) is called \emph{$w$-admissible} if 
\begin{equation} \label{eq:def_admissible}
\Pi(w, A): w_{0} \xrightarrow{|\gamma_{i_{1}}|} w_{1} \xrightarrow{|\gamma_{i_{2}}|} \cdots \xrightarrow{|\gamma_{i_{t}}|} w_{t}
\end{equation}
is a directed path in $\QBG(W)$. 
We denote by $\CA(w, \Gamma)$ the set of all $w$-admissible subsets of $\{1, \ldots, r\}$. 
\end{defn}

We define statistics of admissible subsets used in the Chevalley formula. 
Let $\lambda \in P$ and $w \in W$. Take a $\lambda$-chain $\Gamma = (\gamma_{1}, \ldots, \gamma_{r})$ with corresponding alcove path $(A_{0}, \ldots, A_{r})$. 
Fix a $w$-admissible subset $A = \{i_{1} < \cdots < i_{t}\} \subset \{1, \ldots, r\}$. 
We set 
\begin{equation}
n(A) := \# \{j \in A \mid \gamma_{j} \in -\Delta^{+}\}. 
\end{equation}
Take the corresponding directed path $\Pi(w, A)$ in $\QBG(W)$ defined in \eqref{eq:def_admissible}. 
We define $\ed(A) \in W$ and $\down(A) \in Q^{\vee}$ by 
\begin{align}
\ed(A) &:= w_{t}, \\
\down(A) &:= \sum_{\substack{1 \le k \le t \\ \text{$w_{k-1} \rightarrow w_{k}$ is a quantum edge}}} |\gamma_{i_{k}}|^{\vee}. 
\end{align}
Then we define the \emph{classical part} $\CA|_{Q=0}(w, \Gamma)$ of $\CA(w, \Gamma)$ as 
\begin{equation}
\CA|_{Q=0}(w, \Gamma) := \{A \in \CA(w, \Gamma) \mid \down(A) = 0\}. 
\end{equation}
This is just a combinatorial object for the theory of the (classical) alcove model introduced by Lenart and Postnikov (\cite{LP}). 

In addition, we need the statistic $\wt(A)$. 
For $\alpha \in \Delta$ and $k \in \BZ$, we denote by $s_{\alpha, k}$ the reflection with respect to $H_{\alpha, k}$, defined as $s_{\alpha, k}(\xi) := \xi - (\pair{\xi}{\alpha^{\vee}} - k) \alpha$ for $\xi \in \Fh_{\BR}^{\ast}$. 
For $1 \le k \le r$, take $l_{k} \in \BZ$ such that the common wall of $A_{k-1}$ and $A_{k}$ is contained in $H_{\gamma_{k}, -l_{k}}$. 
Set $\hr_{k} := s_{\gamma_{k}, -l_{k}}$ for $1 \le k \le r$. 
Then we define $\wt(A) \in P$ as 
\begin{equation}
\wt(A) := -w \hr_{i_{1}} \cdots \hr_{i_{t}} (-\lambda). 
\end{equation}

\subsection{The Chevalley formula and line bundle correspondence}

The expansion identity of a class of a semi-infinite Schubert class twisted by a class of a line bundle into an infinite linear combination of semi-infinite Schubert classes is called the \emph{Chevalley formula}. 
First, we recall the combinatorial description of the Chevalley formula by Lenart, Naito, and Sagaki (\cite{LNS}). 
A sequence $\chi = (\chi_{1}, \ldots, \chi_{l})$ of integers is called a \emph{partition} if $\chi_{1} \ge \cdots \ge \chi_{l} > 0$. 
In this case, we call $l$ the \emph{length} of $\chi$, denoted by $\ell(\chi)$. 
Also, we allow the empty partition $\chi = \emptyset$ with $\ell(\chi) = 0$. 
In this case, we understand $\chi = (0)$. 
For $\lambda = \sum_{i \in I} m_{i} \vpi_{i} \in P$, we define $\Par(\lambda)$ by 
\begin{equation}
\Par(\lambda) := \left\{ \bchi = (\bchi_{i})_{i \in I} \ \middle| \ \text{$\bchi_{i}$ is a partition with $\ell(\bchi_{i}) \le \max\{m_{i}, 0\}$ for $i \in I$}\right\}. 
\end{equation}
For $\bchi = (\bchi_{i})_{i \in I} \in \Par(\lambda)$ with $\bchi_{i} = (\bchi_{i}^{(1)} \ge \cdots \ge \bchi_{i}^{(l_{i})})$ (note that $l_{i} = \ell(\bchi_{i})$), we set 
\begin{equation}
\iota(\bchi) := \sum_{i \in I} \bchi_{i}^{(1)} \alpha_{i}^{\vee} \in Q^{\vee, +}; 
\end{equation}
here if $\bchi_{i} = \emptyset$ is an empty partition, then we understand that $\bchi_{i}^{(1)} = 0$. 
Then a combinatorial description of the Chevalley formula is as follows. 

\begin{thm}[{\cite[Theorem~33]{LNS}}]
Let $\lambda \in P$. Take a reduced $\lambda$-chain $\Gamma$. 
Also, let $x = wt_{\xi} \in W_{\af}$. 
In $K_{T}(\bQG)$, we have 
\begin{align}
[\CO_{\bQG(x)}(-w_{\circ}\lambda)] &= \sum_{A \in \CA(w, \Gamma)} \sum_{\bchi \in \Par(\lambda)} (-1)^{n(A)} \e^{\wt(A)} [\CO_{\bQG(\ed(A)t_{\xi + \down(A) + \iota(\bchi)})}] \\ 
&= \sum_{\bchi \in \Par(\lambda)} \st_{\iota(\bchi)} \sum_{A \in \CA(w, \Gamma)} (-1)^{n(A)} \e^{\wt(A)} [\CO_{\bQG(\ed(A)t_{\xi + \down(A)})}]
\end{align}
\end{thm}

We hereby focus on the case $G = \Sp_{2n}(\BC)$. Note that $w_{\circ}\lambda = -\lambda$ for $\lambda \in P$ in this case. 
We concentrate on the case $\lambda = \pm \ve_{J}$ for some $J \subset [1, n]$. 

If $\lambda = \ve_{J}$ for $J \subset [1, n]$, then we have 
\begin{equation}
\Par(\ve_{J}) = \left\{ (\bchi_{j})_{1 \le j \le n} \ \middle| \ \parbox{19em}{$\bchi_{j}$ is a partition with $\ell(\bchi_{j}) \le 1$ for $j \in L_{J}$, and $\bchi_{j} = \emptyset$ for $j \notin L_{J}$} \right\} \simeq \BZ_{\ge 0}^{|L_{J}|} 
\end{equation}
and for $\bchi = (\bchi_{j})_{1 \le j \le n} \in \Par(\ve_{J})$ with $\bchi_{j} = (m_{j})$ and $m_{j} \in \BZ_{\ge 0}$ for $j \in L_{J}$, we have 
\begin{equation}
\iota(\bchi) = \sum_{j \in L_{J}} m_{j} \alpha_{j}^{\vee}.
\end{equation} 

Let us consider the case $\lambda = -\ve_{J}$ for $J \subset [1, n]$. Then we have 
\begin{equation}
\Par(-\ve_{J}) = \left\{ (\bchi_{j})_{1 \le j \le n} \ \middle| \ \parbox{19em}{$\bchi_{j}$ is a partition with $\ell(\bchi_{j}) \le 1$ for $j \in M_{J}$, and $\bchi_{j} = \emptyset$ for $j \notin M_{J}$} \right\} \simeq \BZ_{\ge 0}^{|M_{J}|} 
\end{equation}
and for $\bchi = (\bchi_{j})_{1 \le j \le n} \in \Par(-\ve_{J})$ with $\bchi_{j} = (m_{j})$ and $m_{j} \in \BZ_{\ge 0}$ for $j \in M_{J}$, we have 
\begin{equation}
\iota(\bchi) = \sum_{j \in M_{J}} m_{j} \alpha_{j}^{\vee}. 
\end{equation}

To finish computing the Chevalley formula, we study the structure of certain admissible subsets. 
We use the following fact, which will be proved in Appendix~\ref{sec:appendix}. 
\begin{prop} \label{prop:QAM=AM}
Assume that $G = \Sp_{2n}(\BC)$. Take $J \subset [1, n]$. 
Set $\lambda = \ve_{J}$ or $\lambda = -\ve_{J}$ and take a reduced $\lambda$-chain $\Gamma$. 
Then we have 
\begin{equation}
\sum_{A \in \CA(e, \Gamma)} (-1)^{n(A)} \e^{\wt(A)} [\CO_{\bQG(\ed(A)t_{\down(A)})}] = \sum_{A \in \CA|_{Q=0}(e, \Gamma)} (-1)^{n(A)} \e^{\wt(A)} [\CO_{\bQG(\ed(A))}]. 
\end{equation}
\end{prop}

In the following corollary, we write $1/(1-\st_{j}) := \sum_{m=0}^{\infty} \st_{j}^{m}$ for $1 \le j \le n$, which is a well-defined operator on $K_{T}(\bQG)$.   

\begin{cor}
Assume that $G = \Sp_{2n}(\BC)$. Take $J \subset [1, n]$. 
\begin{enu}
\item For a reduced $\ve_{J}$-chain $\Gamma$, we have 
\begin{align}
[\CO_{\bQG}(\ve_{J})] &= \left( \prod_{j \in L_{J}} \sum_{m_{j}=0}^{\infty} \st_{m_{j}\alpha_{j}^{\vee}} \right) \sum_{A \in \CA|_{Q=0}(e, \Gamma)} (-1)^{n(A)} \e^{\wt(A)} [\CO_{\bQG(\ed(A))}] \\ 
&= \left( \prod_{j \in L_{J}} \frac{1}{1-\st_{j}} \right) \sum_{A \in \CA|_{Q=0}(e, \Gamma)} (-1)^{n(A)} \e^{\wt(A)} [\CO_{\bQG(\ed(A))}]. \label{eq:veJ}
\end{align}

\item For a reduced $(-\ve_{J})$-chain $\Gamma$, we have 
\begin{align}
[\CO_{\bQG}(-\ve_{J})] &= \left( \prod_{j \in M_{J}} \sum_{m_{j}=0}^{\infty} \st_{m_{j}\alpha_{j}^{\vee}} \right) \sum_{A \in \CA|_{Q=0}(e, \Gamma)} (-1)^{n(A)} \e^{\wt(A)} [\CO_{\bQG(\ed(A))}] \\ 
&= \left( \prod_{j \in M_{J}} \frac{1}{1-\st_{j}} \right) \sum_{A \in \CA|_{Q=0}(e, \Gamma)} (-1)^{n(A)} \e^{\wt(A)} [\CO_{\bQG(\ed(A))}]. \label{eq:-veJ}
\end{align}
\end{enu}
\end{cor}

On the other hand, Lenart and Postnikov described the Chevalley formula for the (classical) $K$-ring $K_{T}(G/B)$ for an arbitrary $G$ as follows. 
\begin{prop}[{\cite[Theorem~6.1]{LP}}]
Let $\lambda \in P$ and take a reduced $\lambda$-chain $\Gamma$. Also, let $w \in W$. In $K_{T}(G/B)$, we have 
\begin{equation} \label{eq:classical_Chevalley}
[\CO_{G/B}(\lambda)] \cdot [\CO^{w}] = \sum_{A \in \CA|_{Q=0}(w, \Gamma)} (-1)^{n(A)} \e^{-\wt(A)} [\CO^{\ed(A)}]. 
\end{equation}
\end{prop}
\begin{rem}
In \cite{LP}, the Chevalley formula is described for classes $[\CO_{X_{w}}]$ of Schubert varieties $X_{w} = \overline{B \dot{w} B/B}$, $w \in W$, instead of opposite Schubert classes $[\CO^{w}]$. By using the relation $X^{w} = w_{\circ} X_{w_{\circ}w}$ for $w \in W$, we obtain \eqref{eq:classical_Chevalley}. 

\end{rem}

By combining \eqref{eq:veJ}, \eqref{eq:-veJ}, and \eqref{eq:classical_Chevalley}, we obtain the following correspondence. 

\begin{prop}
Assume that $G = \Sp_{2n}(\BC)$. For $J \subset [1, n]$ and $\CZ \in QK_{T}(\FlC{n})$, we have 
\begin{align}
\Phi([\CO_{G/B}(-\ve_{J})] \star \CZ) &= \left( \prod_{\substack{1 \le j \le n-1 \\ j \notin J, \, j+1 \in J}} (1-\st_{j}) \right) [\CO_{\bQG}(-\ve_{I})] \otimes \Phi(\CZ), \label{eq:correspondence_veJ}\\ 
\Phi([\CO_{G/B}(\ve_{J})] \star \CZ) &= \left( \prod_{\substack{2 \le j \le n \\ j \notin J, \, j-1 \in J}} (1-\st_{j-1}) \right) [\CO_{\bQG}(\ve_{J})] \otimes \Phi(\CZ). \label{eq:correspondence_-veJ}
\end{align}
\end{prop}

\begin{proof}
We show only \eqref{eq:correspondence_veJ} because we can show \eqref{eq:correspondence_-veJ} by the same argument. 
Since $QK_{T}(\FlC{n})$ is generated by line bundle classes $[\CO_{G/B}(-\vpi_{j})]$ for $1 \le j \le n$ as $R(T)\bra{Q_{1}, \ldots, Q_{n}}$-algebra (see \cite[\S 5]{MNS}), it is sufficient to show that 
\begin{equation}
\Phi([\CO_{G/B}(-\ve_{J})]) = \left( \prod_{\substack{1 \le j \le n-1 \\ j \notin I, \, j+1 \in I}} (1-\st_{j}) \right) [\CO_{\bQG}(-\ve_{J})] 
\end{equation}
(cf. Theorem~\ref{thm:QK_vs_semi-infinite}\,(ii)). 
Let $\Gamma$ be a reduced $(-\ve_{J})$-chain. By \eqref{eq:classical_Chevalley}, we see that 
\begin{equation}
[\CO_{G/B}(-\ve_{J})] = [\CO_{G/B}(-\ve_{J})] \cdot [\CO^{e}] = \sum_{A \in \CA|_{Q=0}(e, \Gamma)} (-1)^{n(A)} \e^{-\wt(A)} [\CO^{\ed(A)}]. 
\end{equation}
in $K_{T}(\FlC{n}) \subset QK_{T}(\FlC{n})$. By applying $\Phi$ in Theorem~\ref{thm:QK_vs_semi-infinite}, we obtain 
\begin{align}
\Phi([\CO_{G/B}(-\ve_{J})]) &= \sum_{A \in \CA|_{Q=0}(e, \Gamma)} (-1)^{n(A)} \e^{\wt(A)} [\CO_{\bQG(\ed(A))}] \\ 
&= \left(\prod_{j \in M_{J}} (1-\st_{j}) \right) [\CO_{\bQG}(-\ve_{J})] \qquad \text{by \eqref{eq:veJ}} \\ 
&= \left( \prod_{\substack{1 \le j \le n-1 \\ j \notin I, \, j+1 \in I}} (1-\st_{j}) \right) [\CO_{\bQG}(-\ve_{J})], 
\end{align}
as desired. 
\end{proof}

By using the correspondence \eqref{eq:correspondence_veJ} and \eqref{eq:correspondence_-veJ}, we obtain the following identities, which will be used in the computation of $QK$-Whitney relations. 
\begin{cor} \label{cor:multiple_line_bundles}
Assume that $G = \Sp_{2n}(\BC)$. 
Let $2 \le k \le n$ and take $J \subset [1, k-1]$. Then we have 
\begin{align}
[\CO_{G/B}(-\ve_{J \sqcup \{k\}})] &= \begin{cases} \dfrac{1}{1-Q_{k-1}} [\CO_{G/B}(-\ve_{J})] \star [\CO_{G/B}(-\ve_{k})] & \text{if $k-1 \in J$}, \\[5mm]
[\CO_{G/B}(-\ve_{J})] \star [\CO_{G/B}(-\ve_{k})] & \text{if $k-1 \notin J$}, \end{cases} \label{eq:multiple_1} \\ 
[\CO_{G/B}(\ve_{J \sqcup \{k\}})] &= \begin{cases} \dfrac{1}{1-Q_{k-1}} [\CO_{G/B}(\ve_{J})] \star [\CO_{G/B}(\ve_{k})] & \text{if $k-1 \in J$}, \\[5mm]
[\CO_{G/B}(\ve_{J})] \star [\CO_{G/B}(\ve_{k})] & \text{if $k-1 \notin J$}. \end{cases} \label{eq:multiple_2}
\end{align}
\end{cor}
\begin{proof}
We show only \eqref{eq:multiple_1} because we can show \eqref{eq:multiple_2} by the same argument. 
By \eqref{eq:correspondence_veJ}, we compute 
\begin{align}
& \Phi([\CO_{G/B}(-\ve_{J})] \star [\CO_{G/B}(-\ve_{k})]) \\ 
&= \left( \left( \prod_{\substack{1 \le j \le n-1 \\ j \notin J, j+1 \in J}} (1-\st_{j}) \right) [\CO_{\bQG}(-\ve_{J})] \right) \otimes \left( (1-\st_{k-1}) [\CO_{\bQG}(-\ve_{k})] \right) \\ 
&= \begin{cases}
(1-\st_{k-1}) \left( \displaystyle \prod_{\substack{1 \le j \le n-1 \\ j \notin J \sqcup \{k\} , j+1 \in J \sqcup \{k\}}} (1-\st_{j}) \right) [\CO_{\bQG}(-\ve_{J \sqcup \{k\}})] & \text{if $k-1 \in J$}, \\ 
\left( \displaystyle \prod_{\substack{1 \le j \le n-1 \\ j \notin J \sqcup \{k\} , j+1 \in J \sqcup \{k\}}} (1-\st_{j}) \right) [\CO_{\bQG}(-\ve_{J \sqcup \{k\}})] & \text{if $k-1 \notin J$}. 
\end{cases} \\ 
&= \begin{cases}
\Phi((1-Q_{k-1}) [\CO_{G/B}(-\ve_{J \sqcup \{k\}})]) & \text{if $k-1 \in J$}, \\ 
\Phi([\CO_{G/B}(-\ve_{J \sqcup \{k\}})]) & \text{if $k-1 \notin J$}. 
\end{cases}
\end{align}
Since $\Phi$ is an isomorphism of $R(T)$-modules, we obtain the desired identity \eqref{eq:multiple_1}. 
\end{proof}

\section{The proof of quantum \texorpdfstring{$K$}{K}-theoretic Whitney relations} \label{sec:QK-Whitney_computation}

\subsection{Proof of Theorem~\ref{thm:QK-Whitney_A}}
First, we give a proof \eqref{eq:QK-Whitney_1}. 
For this purpose, we show the following identity. 
\begin{lem}
In $QK_{T}(\FlC{n})$, for $1 \le k \le n$ and $1 \le d \le k$, we have 
\begin{equation} \label{eq:QK-Whitney_lem_1}
\left[ \extprod^{d} \CS_{k} \right] = \left[ \extprod^{d} \CS_{k-1} \right] + \frac{1}{1-Q_{k-1}} [\CS_{k}/\CS_{k-1}] \star \left( \left[ \extprod^{d-1} \CS_{k-1} \right] - Q_{k-1} \left[ \extprod^{d-1} \CS_{k-2} \right] \right). 
\end{equation}
\end{lem}
\begin{proof}
We compute as follows: 
\begin{align}
& \left[ \extprod^{d} \CS_{k} \right] \\ 
& = \sum_{\substack{J \subset [1, k] \\ |J| = d}} [\CO_{G/B}(-\ve_{J})] \\ 
& = \sum_{\substack{J \subset [1, k-1] \\ |J| = d}} [\CO_{G/B}(-\ve_{J})] + \sum_{\substack{J \subset [1, k-1] \\ |J| = d-1}} [\CO_{G/B}(-\ve_{J \sqcup \{k\}})] \\ 
\begin{split}
& = \left[ \extprod^{d} \CS_{k-1} \right] + \sum_{\substack{J \subset [1, k-2] \\ |J| = d-1}} [\CO_{G/B}(-\ve_{k})] \star [\CO_{G/B}(-\ve_{J})] \\ 
& \qquad + \frac{1}{1-Q_{k-1}} \sum_{\substack{J \subset [1, k-1] \\ k-1 \in J \\ |J| = d-1}} [\CO_{G/B}(-\ve_{k})] \star [\CO_{G/B}(-\ve_{J})]
\end{split} \qquad \text{by \eqref{eq:multiple_1}} \\ 
\begin{split}
& = \left[ \extprod^{d} \CS_{k-1} \right] + \frac{1}{1-Q_{k-1}} [\CO_{G/B}(-\ve_{k})] \star \sum_{\substack{J \subset [1, k-1] \\ |J| = d-1}} [\CO_{G/B}(-\ve_{J})] \\ 
& \qquad - \frac{Q_{k-1}}{1-Q_{k-1}} [\CO_{G/B}(-\ve_{k})] \star \sum_{\substack{J \subset [1, k-2] \\ |J| = d-1}} [\CO_{G/B}(-\ve_{J})]
\end{split} \\ 
\begin{split}
& = \left[ \extprod^{d} \CS_{k-1} \right] \\ 
& \qquad + \frac{1}{1-Q_{k-1}} [\CO_{G/B}(-\ve_{k})] \star \left( \sum_{\substack{J \subset [1, k-1] \\ |J| = d-1}} [\CO_{G/B}(-\ve_{J})] - Q_{k-1} \sum_{\substack{J \subset [1, k-2] \\ |J| = d-1}} [\CO_{G/B}(-\ve_{J})] \right)
\end{split} \\ 
&= \left[ \extprod^{d} \CS_{k-1} \right] + \frac{1}{1-Q_{k-1}} [\CS_{k}/\CS_{k-1}] \star \left( \left[ \extprod^{d-1} \CS_{k-1} \right] - Q_{k-1} \left[ \extprod^{d-1} \CS_{k-2} \right] \right), 
\end{align}
as desired. 
\end{proof}

\begin{proof}[Proof of Theorem~\ref{thm:QK-Whitney_A}\,(1)]
The left-hand side of \eqref{eq:QK-Whitney_1} is computed as follows:  
\begin{align}
& \lambda_{y}(\CS_{k-1}) \star \lambda_{y}(\CS_{k}/\CS_{k-1}) \\ 
&= \left( \sum_{d = 0}^{k-1} y^{d} \left[ \extprod^{d} \CS_{k-1} \right] \right) \star (1 + y [\CS_{k}/\CS_{k-1}]) \\ 
&= \sum_{d = 0}^{k} y^{d} \left( \left[ \extprod^{d} \CS_{k-1} \right] + \left[ \extprod^{d-1} \CS_{k-1} \right] \star [\CS_{k}/\CS_{k-1}] \right). \label{eq:computation_LHS_QK-Whitney_1}
\end{align}

On the other hand, the right-hand side of \eqref{eq:QK-Whitney_1} is computed as follows: 
\begin{align}
& \lambda_{y}(\CS_{k}) - y \frac{Q_{k-1}}{1-Q_{k-1}} [\CS_{k}/\CS_{k-1}] \star (\lambda_{y}(\CS_{k-1}) - \lambda_{y}(\CS_{k-2})) \\ 
&= \sum_{d=0}^{k} y^{d} \left[ \extprod^{d} \CS_{k} \right] - y \frac{Q_{k-1}}{1-Q_{k-1}} [\CS_{k}/\CS_{k-1}] \star \left( \sum_{d=0}^{k-1} y^{d} \left[ \extprod^{d} \CS_{k-1} \right] - \sum_{d=0}^{k-2} y^{d} \left[ \extprod^{d} \CS_{k-2} \right] \right) \\ 
&= \sum_{d=0}^{k} y^{d} \left[ \extprod^{d} \CS_{k} \right] - \frac{Q_{k-1}}{1-Q_{k-1}} [\CS_{k}/\CS_{k-1}] \star \left( \sum_{d=1}^{k} y^{d} \left[ \extprod^{d-1} \CS_{k-1} \right] - \sum_{d=1}^{k-1} y^{d} \left[ \extprod^{d-1} \CS_{k-2} \right] \right) \\ 
&= \sum_{d=0}^{k} y^{d} \left( \left[ \extprod^{d} \CS_{k} \right] - \frac{Q_{k-1}}{1-Q_{k-1}} [\CS_{k}/\CS_{k-1}] \star \left( \left[ \extprod^{d-1} \CS_{k-1} \right] - \left[ \extprod^{d-1} \CS_{k-2} \right] \right) \right) \label{eq:computation_RHS_QK-Whitney_1}
\end{align}
By \eqref{eq:QK-Whitney_lem_1}, we conclude that $\eqref{eq:computation_LHS_QK-Whitney_1} = \eqref{eq:computation_RHS_QK-Whitney_1}$, as desired. 
\end{proof}

Then \eqref{eq:QK-Whitney_2} can be proved in parallel with \eqref{eq:QK-Whitney_1} by showing the following identity analogous to \eqref{eq:QK-Whitney_lem_1}, of which we omit the proof. 

\begin{lem}
In $QK_{T}(\FlC{n})$, for $1 \le d \le k \le n$ and $1 \le d \le k$, we have 
\begin{equation} \label{eq:QK-Whitney_lem_2}
\left[ \extprod^{d} \CS_{k}^{\vee} \right] = \left[ \extprod^{d} \CS_{k-1}^{\vee} \right] + \frac{1}{1-Q_{k-1}} [(\CS_{k}/\CS_{k-1})^{\vee}] \star \left( \left[ \extprod^{d-1} \CS_{k-1}^{\vee} \right] - Q_{k-1} \left[ \extprod^{d-1} \CS_{k-2}^{\vee} \right] \right). 
\end{equation}
\end{lem}

\subsection{The Borel-type relation}
To show Theorem~\ref{thm:QK-Whitney_B}, we need to review the Borel-type relations of $QK_{T}(\FlC{n})$ described in \cite{KN2}, which are enough to construct a defining ideal of $QK_{T}(\FlC{n})$; that is, an ideal $I_{Q}^{\SB}$ of a Laurent polynomial ring $(R(T)\bra{Q_{1}, \ldots, Q_{n}})[z_{1}^{\pm 1}, \ldots, z_{n}^{\pm 1}]$ such that
\begin{equation}
QK_{T}(\FlC{n}) \simeq (R(T)\bra{Q_{1}, \ldots, Q_{n}})[z_{1}^{\pm 1}, \ldots, z_{n}^{\pm 1}]/I_{Q}^{\SB}. 
\end{equation}
as $R(T)$-algebras. This isomorphism is called the \emph{Borel-type presentation}. 

\begin{defn}[{\cite[Definition~6.1]{KN2}}]
Let $J \subset [1, \overline{1}]$. 
\begin{enu}
\item For $1 \le j \le n$, we define $\varphi_{J}(j) \in \BZ\bra{Q_{1}, \ldots, Q_{n}}$ by 
\begin{equation}
\varphi_{J}(j) := \begin{cases}
\dfrac{1}{1-Q_{j}} & \text{if $j, j+1 \in J$}, \\ 
0 & \text{otherwise}. 
\end{cases}
\end{equation}

\item For $2 \le j \le n$, we define $\varphi_{J}(\overline{j}) \in \BZ\bra{Q_{1}, \ldots, Q_{n}}$ by 
\begin{equation}
\varphi_{J}(\overline{j}) := \begin{cases}
1 + \dfrac{Q_{j-1}Q_{j} \cdots Q_{n}}{1-Q_{j-1}} & \text{if $J = \{ \cdots < j-1 < \overline{j-1} < \cdots \}$}, \\ 
\dfrac{1}{1-Q_{j-1}} & \text{if $\overline{j}, \overline{j-1} \in J$}, \\ 
1 & \text{otherwise}. 
\end{cases}
\end{equation}

\item We set $\varphi_{J}(\overline{1}) := 1 \in \BZ\bra{Q_{1}, \ldots, Q_{n}}$. 
\end{enu}
\end{defn}

Throughout this paper, we denote by $\prod\nolimits^{\star}$ the product with respect to the quantum product $\star$ on $QK_{T}(\FlC{n})$. 
Then the Borel-type relations are described as follows; here recall that $\ve_{\overline{j}} = -\ve_{j}$ for $1 \le j \le n$. 
\begin{prop}[{\cite[Corollary~6.4]{KN2}}] \label{prop:QK-Borel_rel}
In $QK_{T}(\FlC{n})$, for $0 \le d \le 2n$, we have 
\begin{equation} \label{eq:QK-Borel_rel}
\sum_{\substack{J \subset [1, \overline{1}] \\ |J|=d}} \left( \prod_{1 \le j \le \overline{1}} \varphi_{J}(j) \right) \left( \qprod_{j \in J} [\CO_{G/B}(-\ve_{j})] \right) = \left[ \extprod^{d} \udBC^{2n} \right]. 
\end{equation}
\end{prop}

Also, we use the following relation. 
\begin{prop}[{\cite[Lemma~6.7]{KN2}}]
In $QK_{T}(\FlC{n})$, for $1 \le j \le n$, we have
\begin{equation} \label{eq:quantum_inverse}
[\CO_{G/B}(-\ve_{j})] \star [\CO_{G/B}(\ve_{j})] = (1-Q_{j-1})(1-Q_{j}). 
\end{equation}
\end{prop}

\subsection{Proof of Theorem~\ref{thm:QK-Whitney_B}}
We prove \eqref{eq:QK-Whitney_3} by comparing coefficients of $y^{d}$ of both sides for each $1 \le d \le 2n$. 
\begin{lem}
In $QK_{T}(\FlC{n})$, for $1 \le d \le 2n$, we have
\begin{equation} \label{eq:QK-Whitney_computation_lem1}
\begin{split}
& \sum_{\substack{k, l \ge 0 \\ k+l = d}} \left[ \extprod^{k} \CS_{n} \right] \star \left[ \extprod^{l} \CS_{n}^{\vee} \right] = \left[ \extprod^{d} \udBC^{2n} \right] \\ 
& \hspace{20mm} - \sum_{p = 1}^{n} \frac{Q_{p}Q_{p+1} \cdots Q_{n}}{1 - Q_{p}} \sum_{\substack{J, K \subset [1, n] \\ |J| + |K| = d \\ \max J = \max K = p}} [\CO_{G/B}(-\ve_{J})] \star [\CO_{G/B}(\ve_{K})]. 
\end{split}
\end{equation}
\end{lem}
\begin{proof}
In what follows, we set $\overline{K} := \{\overline{j} \mid j \in K \}$ for $K \subset [1, n]$. 
Let $J, K \subset [1, n]$. 
By \eqref{eq:multiple_1}, we have 
\begin{align}
[\CO_{G/B}(-\ve_{J})] &= \left( \prod_{\substack{2 \le j \le n \\ j-1 \in J, j \in J}} \frac{1}{1-Q_{j-1}} \right) \left( \qprod_{j \in J} [\CO_{G/B}(-\ve_{j})] \right) \\ 
&= \left( \prod_{1 \le j \le n} \varphi_{J}(j) \right) \left( \qprod_{j \in J} [\CO_{G/B}(-\ve_{j})] \right). 
\end{align}
Also, by \eqref{eq:multiple_2}, we have 
\begin{align}
[\CO_{G/B}(\ve_{K})] &= \left( \prod_{\substack{2 \le j \le n \\ j-1\in K, j \in K}} \frac{1}{1-Q_{j-1}} \right) \left( \qprod_{j \in K} [\CO_{G/B}(\ve_{j})] \right) \\ 
&= \left( \prod_{\overline{n} \le j \le \overline{1}} \varphi_{\overline{K}}(j) \right) \left( \qprod_{j \in \overline{K}} [\CO_{G/B}(-\ve_{j})] \right). 
\end{align}
Hence we have 
\begin{equation}
[\CO_{G/B}(-\ve_{J})] \star [\CO_{G/B}(\ve_{K})] = \left( \prod_{1 \le j \le n} \varphi_{J}(j) \right) \left( \prod_{\overline{n} \le j \le \overline{1}} \varphi_{\overline{K}}(j) \right) \left( \ \qprod_{j \in J \sqcup \overline{K}} [\CO_{G/B}(-\ve_{j})] \right). 
\end{equation}

If $\max J \not= \max K$, then we have 
\begin{equation}
\varphi_{J \sqcup \overline{K}}(j) = \begin{cases}
\varphi_{J}(j) & \text{if $1 \le j \le n$}, \\ 
\varphi_{\overline{K}}(j) & \text{if $\overline{n} \le j \le \overline{1}$}. 
\end{cases}
\end{equation}
Hence we obtain 
\begin{equation}
\left( \prod_{1 \le j \le n} \varphi_{J}(j) \right) \left( \prod_{\overline{n} \le j \le \overline{1}} \varphi_{\overline{K}}(j) \right) = \prod_{1 \le j \le \overline{1}} \varphi_{J \sqcup \overline{K}} (j). 
\end{equation}

Assume that $\max J = \max K =: p$. If $p = n$, then we see that 
\begin{equation}
\varphi_{J \sqcup \overline{K}}(j) = \begin{cases}
\varphi_{J}(j) & \text{if $1 \le j \le n-1$}, \\ 
\dfrac{1}{1-Q_{n}} = \left( 1 + \dfrac{Q_{n}}{1-Q_{n}} \right)\underbrace{\varphi_{J}(n)}_{=1} & \text{if $j = n$} \\ 
\varphi_{\overline{K}}(j) & \text{if $\overline{n} \le j \le \overline{1}$}. 
\end{cases}
\end{equation}
If $p < n$, then 
\begin{equation}
\varphi_{J \sqcup \overline{K}} (j) = \begin{cases}
\varphi_{J}(j) & \text{if $1 \le j \le n$}, \\ 
\varphi_{\overline{K}}(j) & \text{if $\overline{n} \le j < \overline{p+1}$ or $\overline{p} \le j \le \overline{1}$}, \\ 
\left( 1 + \dfrac{Q_{p} Q_{p+1} \cdots Q_{n}}{1 - Q_{p}} \right) \underbrace{\varphi_{\overline{K}}(\overline{p+1})}_{=1} & \text{if $j = \overline{p+1}$}. 
\end{cases}
\end{equation}
Hence in the case $\max J = \max K = p \in [1, n]$, we obtain 
\begin{equation}
\begin{split}
& \left( \prod_{1 \le j \le n} \varphi_{J}(j) \right) \left( \prod_{\overline{n} \le j \le \overline{1}} \varphi_{\overline{K}}(j) \right) \\ 
& \quad = \prod_{1 \le j \le \overline{1}} \varphi_{J \sqcup \overline{K}}(j) - \frac{Q_{p}Q_{p+1} \cdots Q_{n}}{1-Q_{p}} \left( \prod_{1 \le j \le n} \varphi_{J}(j) \right) \left( \prod_{\overline{n} \le j \le \overline{1}} \varphi_{\overline{K}}(j) \right). 
\end{split}
\end{equation}
Therefore, we compute 
\begin{align}
& \sum_{\substack{k, l \ge 0 \\ k+l=d}} \left[ \extprod^{k} \CS_{n} \right] \star \left[ \extprod^{l} \CS_{n}^{\vee} \right] \\ 
&= \sum_{\substack{k, l \ge 0 \\ k+l=d}} \sum_{\substack{J, K \subset [1, n] \\ |J|=k, |K|=l}} [\CO_{G/B}(-\ve_{J})] \star [\CO_{G/B}(\ve_{K})] \\ 
\begin{split}
&= \sum_{\substack{L \subset [1, \overline{1}] \\ |L|=d}} \left( \prod_{1 \le j \le \overline{1}} \varphi_{L}(j) \right) \left( \qprod_{j \in L} [\CO_{G/B}(-\ve_{j})] \right) \\ 
& \qquad - \sum_{p=1}^{n} \sum_{\substack{k, l \ge 0 \\ k+l = d}}\sum_{\substack{J, K \subset [1, n] \\ |J|=k, |K|=l, \\ \max J = \max K = p}}
\frac{Q_{p}Q_{p+1} \cdots Q_{n}}{1-Q_{p}} \\ 
& \hspace{40mm} \times \left( \prod_{1 \le j \le n} \varphi_{J}(j) \right) \left( \prod_{\overline{n} \le j \le \overline{1}} \varphi_{\overline{K}}(j) \right) \left( \ \qprod_{j \in J \sqcup \overline{K}} [\CO_{G/B}(-\ve_{j})] \right)
\end{split} \\ 
&= \left[ \extprod^{d} \udBC^{2n} \right] - \sum_{p=1}^{n} \sum_{\substack{J, K \subset [1, n] \\ |J|+|K|=d, \\ \max J = \max K = p}}
\frac{Q_{p}Q_{p+1} \cdots Q_{n}}{1-Q_{p}} [\CO_{G/B}(-\ve_{J})] \star [\CO_{G/B}(\ve_{K})], 
\end{align}
as desired; here in the final equality, we used \eqref{eq:QK-Borel_rel}. 
\end{proof}

\begin{lem}
In $QK_{T}(\FlC{n})$, for $1 \le d \le 2n$ and $1 \le p \le n$, we have
\begin{equation} \label{eq:QK-Whitney_computation_lem2}
\begin{split}
& \sum_{\substack{J, K \subset [1, n] \\ |J| + |K| = d \\ \max J = \max K = p}} [\CO_{G/B}(-\ve_{J})] \star [\CO_{G/B}(\ve_{K})] \\ 
& \qquad = \frac{1-Q_{p}}{1-Q_{p-1}} \sum_{\substack{k, l \ge 0 \\ k+l = d-2}} \left( \left[ \extprod^{k} \CS_{p-1} \right] - Q_{p-1} \left[ \extprod^{k} \CS_{p-2} \right] \right) \\ 
& \hspace{60mm} \star \left( \left[ \extprod^{l} \CS_{p-1}^{\vee} \right] - Q_{p-1} \left[ \extprod^{l} \CS_{p-2}^{\vee} \right] \right). 
\end{split}
\end{equation}
\end{lem}

\begin{proof}
Assume that $p \ge 2$. We compute as 
\begin{align}
& \sum_{\substack{J, K \subset [1, n] \\ |J| + |K| = d \\ \max J = \max K = p}} [\CO_{G/B}(-\ve_{J})] \star [\CO_{G/B}(\ve_{K})] \\ 
\begin{split}
&= \sum_{\substack{J, K \subset [1, n] \\ |J| + |K| = d \\ \max J = \max K = p}} \left( \varphi_{J}(p) [\CO_{G/B}(-\ve_{J \setminus \{p\}})] \star [\CO_{G/B}(-\ve_{p})] \right) \\ 
& \hspace{40mm} \star \left( \varphi_{\overline{K}}(\overline{p}) [\CO_{G/B}(\ve_{K \setminus \{p\}})] \star [\CO_{G/B}(\ve_{p})] \right)
\end{split} & \hspace{-20mm} \text{by \eqref{eq:multiple_1} and \eqref{eq:multiple_2}} & \\ 
\begin{split}
&= \sum_{\substack{J, K \subset [1, n] \\ |J| + |K| = d \\ \max I = \max J = p}} (1-Q_{p-1})(1-Q_{p}) \\ 
& \hspace{40mm} \times \left( \varphi_{J}(p) [\CO_{G/B}(-\ve_{J \setminus \{p\}})] \right) \star \left( \varphi_{\overline{K}}(\overline{p}) [\CO_{G/B}(\ve_{K \setminus \{p\}})] \right)
\end{split} & \text{by \eqref{eq:quantum_inverse}} & \\ 
\begin{split}
&= (1-Q_{p-1})(1-Q_{p}) \\ 
& \hspace{10mm} \times \sum_{\substack{k, l \ge 0 \\ k+l = d-2}} \left( \sum_{\substack{J \subset [1, p-1] \\ p-1 \in J \\ |J|=k}} \frac{1}{1-Q_{p-1}} [\CO_{G/B}(-\ve_{J})] + \sum_{\substack{J \subset [1, p-1] \\ p-1 \notin J \\ |J|=k}} [\CO_{G/B}(-\ve_{J})] \right) \\ 
& \hspace{30mm} \star \left( \sum_{\substack{K \subset [1, p-1] \\ p-1 \in k \\ |K|=l}} \frac{1}{1-Q_{p-1}} [\CO_{G/B}(\ve_{K})] + \sum_{\substack{K \subset [1, p-1] \\ p-1 \notin K \\ |K|=l}} [\CO_{G/B}(\ve_{K})] \right)
\end{split} \\ 
\begin{split}
&= \frac{1-Q_{p}}{1-Q_{p-1}} \sum_{\substack{k, l \ge 0 \\ k+l=d-2}} \left( \sum_{\substack{J \subset [1, p-1] \\ |J|=k}} [\CO_{G/B}(-\ve_{J})] - Q_{p-1} \sum_{\substack{J \subset [1, p-2] \\ |J|=k}} [\CO_{G/B}(-\ve_{J})] \right) \\ 
& \hspace{35mm} \star \left( \sum_{\substack{K \subset [1, p-1] \\ |K|=l}} [\CO_{G/B}(\ve_{K})] - Q_{p-1} \sum_{\substack{K \subset [1, p-2] \\ |K|=l}} [\CO_{G/B}(\ve_{K})] \right)
\end{split} \\ 
&= \frac{1-Q_{p}}{1-Q_{p-1}} \sum_{\substack{k, l \ge 0 \\ k+l=d-2}} \left( \left[ \extprod^{k} \CS_{p-1} \right] - Q_{p-1} \left[ \extprod^{k} \CS_{p-2} \right] \right) \\ 
& \hspace{60mm} \star \left( \left[ \extprod^{l} \CS_{p-1} \right] - Q_{p-1} \left[ \extprod^{l} \CS_{p-2} \right] \right), 
\end{align}
as desired. If $p = 1$, then since
\begin{equation}
[\CO_{G/B}(-\ve_{1})] \star [\CO_{G/B}(\ve_{1})] = 1-Q_{1} 
\end{equation}
by \eqref{eq:quantum_inverse}, we see that 
\begin{equation}
\sum_{\substack{J, K \subset [1, n] \\ |J|+|K|=d \\ \max J = \max K = 1}} [\CO_{G/B}(-\ve_{J})] \star [\CO_{G/B}(\ve_{K})] = \begin{cases}
1-Q_{1} & \text{if $d = 2$}, \\ 
0 & \text{if $d \not= 2$}. 
\end{cases}
\end{equation}
This agrees with \eqref{eq:QK-Whitney_computation_lem2}. 
\end{proof}

By combining \eqref{eq:QK-Whitney_computation_lem1} and \eqref{eq:QK-Whitney_computation_lem2}, we obtain the following. 

\begin{prop} 
In $QK_{T}(\FlC{n})$, for $1 \le d \le 2n$, we have 
\begin{equation} \label{eq:QK-Whitney_lem_3}
\begin{split}
& \sum_{\substack{k, l \ge 0 \\ k+l = d}} \left[ \extprod^{k} \CS_{n} \right] \star \left[ \extprod^{l} \CS_{n}^{\vee} \right] = \left[ \extprod^{d} \udBC^{2n} \right] \\ 
& \qquad - \sum_{p = 1}^{n} \frac{Q_{p} Q_{p+1} \cdots Q_{n}}{1 - Q_{p-1}} \sum_{\substack{k, l \ge 0 \\ k+l = d-2}} \left( \left[ \extprod^{k} \CS_{p-1} \right] - Q_{p-1} \left[ \extprod^{k} \CS_{p-2} \right] \right) \\ 
& \hspace{80mm} \star \left( \left[ \extprod^{l} \CS_{p-1}^{\vee} \right] - Q_{p-1} \left[ \extprod^{l} \CS_{p-2}^{\vee} \right] \right). 
\end{split}
\end{equation}
\end{prop}

\begin{proof}[Proof of Theorem~\ref{thm:QK-Whitney_B}]
The left-hand side of \eqref{eq:QK-Whitney_lem_3} is the coefficient of $y^{d}$ in $\lambda_{y}(\CS_{n}) \star \lambda_{y}(\CS_{n}^{\vee})$. 
On the other hand, the right-hand side of \eqref{eq:QK-Whitney_lem_3} is identical to the coefficient of $y^{d}$ in the right-hand side of \eqref{eq:QK-Whitney_3}; recall that
\begin{equation}
\left[ \extprod^{d} \udBC^{2n} \right] = e_{d}(\e^{\ve_{1}}, \e^{\ve_{2}}, \ldots, \e^{\ve_{n}}, \e^{-\ve_{n}}, \ldots, \e^{-\ve_{2}}, \e^{-\ve_{1}}). 
\end{equation}
Therefore \eqref{eq:QK-Whitney_3} holds. 
\end{proof}

\section{The proof of the Whitney-type presentation} \label{sec:Nakayama}
The purpose of this section is to prove Theorem~\ref{thm:Whitney_presentation}. 
The strategy of the proof is to reduce the Whitney-type presentation to the classical presentation of $K_{T}(\FlC{n})$ via Nakayama-type argument established by Gu, Mihalcea, Sharpe, Xu, Zhang, and Zou (\cite{GMSXZZ2}). 

In this section, we review the Nakayama-type argument, and then we give a proof of our Whitney-type presentation. 

\subsection{Nakayama-type result for the quantum \texorpdfstring{$K$}{K}-ring} 
To prove Theorem~\ref{thm:Whitney_presentation}, we review the argument established by Gu, Mihalcea, Sharpe, Xu, Zhang, and Zou (\cite{GMSXZZ2}), which can be obtained as an application of Nakayama's lemma in the theory of commutative rings. In the following theorem, let $G$ be of arbitrary type. Take a parabolic subgroup $P \subset G$. We denote by $Q_{1}, \ldots, Q_{k}$ the Novikov variables associated to $G/P$. 

\begin{thm}[{\cite[Theorem~4.1]{GMSXZZ2}}] \label{thm:Nakayama}

Assume that we are given an isomorphism of $R(T)$-algebras
\begin{equation}
\Phi: R(T)[m_{1}, \ldots, m_{s}]/\langle P_{1}, \ldots, P_{r} \rangle \xrightarrow{\sim} K_{T}(G/P)
\end{equation}
where $m_{1}, \ldots, m_{s}$ are indeterminates, and the $P_{i}$ are polynomials in the variables $m_{j}$ with coefficients in $R(T)$. 
Consider any power series
\begin{equation}
\ti{P}_{i} = \ti{P}_{i}(m_{1}, \ldots, m_{s}; Q_{1}, \ldots, Q_{k}) \in (R(T)[m_{1}, \ldots, m_{s}])\bra{Q_{1}, \ldots, Q_{k}}
\end{equation}
such that 
\begin{itemize}
\item $\ti{P}_{i}(m_{1}, \ldots, m_{s}; 0, \ldots, 0) = P_{i}(m_{1}, \ldots, m_{s})$ for all $m_{1}, \ldots, m_{s}$ and all $i$, and 
\item $\ti{P}_{i}(m_{1}, \ldots, m_{s}; Q_{1}, \ldots, Q_{k}) = 0$ in $QK_{T}(G/P)$, where the $m_{i}$ are regarded as elements in $QK_{T}(G/P)$ via the isomorphism $\Phi$. 
\end{itemize}
Then $\Phi$ lifts to an isomorphism of $R(T)\bra{Q_{1}, \ldots, Q_{k}}$-algebras
\begin{equation}
\ti{\Phi}: (R(T)\bra{Q_{1}, \ldots, Q_{k}})[m_{1}, \ldots, m_{s}]/\langle \ti{P}_{1}, \ldots, \ti{P}_{r} \rangle \xrightarrow{\sim} QK_{T}(G/P). 
\end{equation}
\end{thm}

Due to Theorem~\ref{thm:Nakayama}, we can reduce an isomorphism for $QK_{T}(\FlC{n})$ to that for $K_{T}(\FlC{n})$. 

\subsection{Proof of Theorem~\ref{thm:Whitney_presentation}}
Assume again that $G = \Sp_{2n}(\BC)$. 
We give a proof of Theorem~\ref{thm:Whitney_presentation}. 
To apply Theorem~\ref{thm:Nakayama}, we should compute restrictions of generators of the ideal $I_{Q}^{\SW}$. 
By substituting $Q_{1} = \cdots = Q_{n} = 0$, we see that 
\eqref{eq:QK-Whitney_def1} goes to 
\begin{equation} \label{eq:Whitney_1}
F_{d}^{k} - (F_{d}^{k-1} + F_{d-1}^{k-1}x_{k}) \quad \text{for $1 \le k \le n$ and $1 \le d \le k$}, 
\end{equation}
\eqref{eq:QK-Whitney_def2} goes to 
\begin{equation} \label{eq:Whitney_2} 
G_{d}^{k} - (G_{d}^{k-1} + G_{d-1}^{k-1}y_{k}) \quad \text{for $1 \le k \le n$ and $1 \le d \le k$}, 
\end{equation}
\eqref{eq:QK-Whitney_def3} goes to 
\begin{equation} \label{eq:Whitney_3}
\sum_{k=0}^{d} F_{k}^{n}G_{d-k}^{n} - e_{d}(\e^{\ve_{1}}, \e^{\ve_{2}}, \ldots, \e^{\ve_{n}}, \e^{-\ve_{n}}, \ldots, \e^{-\ve_{2}}, \e^{-\ve_{1}}) \quad \text{for $1 \le d \le n$}, 
\end{equation}
and \eqref{eq:inverse} goes to 
\begin{equation} \label{eq:Whitney_4}
x_{j}y_{j} = 1 \quad \text{for $1 \le j \le n$}, 
\end{equation}
while \eqref{eq:initial_1} and \eqref{eq:initial_2} are not be changed. 

Let $I^{\SW}$ be an ideal of $R^{\SW} := R(T)[F_{d}^{k}, G_{d}^{k}, x_{j}, y_{j} \mid 1 \le k \le n, \ 1 \le d \le k, \ 1 \le j \le n]$ generated by \eqref{eq:Whitney_1}, \eqref{eq:Whitney_2}, \eqref{eq:Whitney_3}, \eqref{eq:Whitney_4}, and \eqref{eq:initial_1}, \eqref{eq:initial_2}. 
If there is no risk of confusion, we use the same notation $a$ for the residue class $a + I^{\SW}$ for $a \in R^{\SW}$. 
Then in $R^{\SW}/I^{\SW}$, we have $x_{j}y_{j} = 1$, $1 \le j \le n$, by \eqref{eq:Whitney_4}. 
Therefore, we see that $x_{j}$ and $y_{j}$ are invertible and $y_{j} = x_{j}^{-1}$. 
Also, by \eqref{eq:Whitney_1} and \eqref{eq:Whitney_2}, in $R^{\SW}/I^{\SW}$, we have 
\begin{align}
F_{d}^{k} &= F_{d}^{k-1} + F_{d-1}^{k-1}x_{k}, \label{eq:Whitney_1_classical} \\ 
G_{d}^{k} &= G_{d}^{k-1} + G_{d-1}^{k-1}x_{k}^{-1} \label{eq:Whitney_2_classical}
\end{align}
for $1 \le k \le n$ and $1 \le d \le k$. In addition, we have $F_{1}^{1} = x_{1}$ and $G_{1}^{1} = x^{-1}$ by \eqref{eq:initial_1} and \eqref{eq:initial_2}. 

\begin{prop} \label{prop:FG_to_el-sym}
In $R^{\SW}/I^{\SW}$, we have 
\begin{equation} \label{eq:Fd_to_el-sym}
F_{d}^{k} = e_{d}(x_{1}, \ldots, x_{k})
\end{equation}
and 
\begin{equation} \label{eq:Gd_to_el-sym}
G_{d}^{k} = e_{d}(x_{1}^{-1}, \ldots, x_{k}^{-1})
\end{equation}
for $1 \le k \le n$ and $1 \le d \le k$. 
\end{prop}
\begin{proof}
For the proof of $F_{d}^{k} = e_{d}(x_{1}, \ldots, x_{k})$, it is enough to verify that $e_{d}(x_{1}, \ldots, x_{k})$ satisfies the same recursion identity as \eqref{eq:Whitney_1_classical}; that is, 
\begin{equation}
e_{d}(x_{1}, \ldots, x_{k}) = e_{d}(x_{1}, \ldots, x_{k-1}) + e_{d-1}(x_{1}, \ldots, x_{k-1})x_{k}; 
\end{equation}
which is obvious by the definition of elementary symmetric polynomials. 
Since $e_{1}(x_{1}) = x_{1}$, we obtain \eqref{eq:Fd_to_el-sym}. 
By the same argument with \eqref{eq:Whitney_2_classical} and $e_{1}(x_{1}^{-1}) = x_{1}^{-1}$, we obtain \eqref{eq:Gd_to_el-sym}. 
\end{proof}

\begin{prop} \label{prop:FtimesG}
In $R^{\SW}/I^{\SW}$, we have 
\begin{equation} \label{eq:FtimesG}
\sum_{k=0}^{d} F_{k}^{n} G_{d-k}^{n} = e_{d}(x_{1}, x_{2}, \ldots, x_{n}, x_{n}^{-1}, \ldots, x_{2}^{-1}, x_{1}^{-1})
\end{equation}
for $1 \le d \le n$. 
\end{prop}
\begin{proof}
By the definition of elementary symmetric polynomials, we have 
\begin{equation}
\sum_{\substack{k. l \ge 0 \\ k+l = d}} e_{k}(x_{1}, \ldots, x_{n}) e_{l}(x_{1}^{-1}, \ldots, x_{n}^{-1}) = e_{d}(x_{1}, x_{2}, \ldots, x_{n}, x_{n}^{-1}, \ldots, x_{2}^{-1}, x_{1}^{-1}). 
\end{equation}
Therefore, Proposition~\ref{prop:FG_to_el-sym} shows \eqref{eq:FtimesG}.  
\end{proof}

Set $R^{\SB} := R(T)[z_{1}^{\pm 1}, \ldots, z_{n}^{\pm 1}]$ and take an ideal $I^{\SB}$ of $R^{\SB}$ generated by 
\begin{equation}
e_{d}(z_{1}, z_{2}, \ldots, z_{n}, z_{n}^{-1}, \ldots, z_{2}^{-1}, z_{1}^{-1}) - e_{d}(\e^{\ve_{1}}, \e^{\ve_{2}}, \ldots, \e^{\ve_{n}}, \e^{-\ve_{n}}, \ldots, \e^{-\ve_{2}}, \e^{-\ve_{1}}) \quad \text{for $1 \le d \le n$}. 
\end{equation}
By Propositions~\ref{prop:FG_to_el-sym} and \ref{prop:FtimesG}, we deduce the following. 
\begin{cor}
We have an isomorphism
\begin{equation}
\Xi_{0}: R^{\SW}/I^{\SW} \xrightarrow{\sim} R^{\SB}/I^{\SB}
\end{equation}
of $R(T)$-algebras such that $\Xi_{0}(x_{j}) = z_{j}$ for $1 \le j \le n$. 
\end{cor}

\begin{proof}[Proof of Theorem~\ref{thm:Whitney_presentation}]
There exists a well-known isomorphism of $R(T)$-algebras  
\begin{equation}
K_{T}(G/B) \simeq R(T) \otimes_{R(T)^{W}} R(T)
\end{equation}
for $G$ of arbitrary type (cf. \cite{PR}); where $R(T)^{W}$ denotes the ring of $W$-invariant elements of $R(T)$.  
In the case $G = \Sp_{2n}(\BC)$, \eqref{eq:rep_ring} implies that 
\begin{equation}
R^{\SB}/I^{\SB} \simeq R(T) \otimes_{R(T)^{W}} R(T)
\end{equation}
and we obtain an isomorphism of $R(T)$-algebras 
\begin{equation}
\Psi: R^{\SB}/I^{\SB} \xrightarrow{\sim} K_{T}(\FlC{n})
\end{equation}
such that  $\Psi(z_{j}) = [\CO_{G/B}(-\ve_{j})]$ for $1 \le j \le n$. 
Therefore, we obtain an isomorphism
\begin{equation}
\Xi = \Psi \circ \Xi_{0}: R^{\SW}/I^{\SW} \xrightarrow{\sim} K_{T}(\FlC{n})
\end{equation}
of $R(T)$-algebras such that $\Xi(x_{j}) = [\CO_{G/B}(-\ve_{j})]$ for $1 \le j \le n$. 
Then applying Theorem~\ref{thm:Nakayama} to this $\Xi$ completes the proof of the theorem. 
\end{proof}

\appendix

\section{Proof of Proposition~\ref{prop:QAM=AM}} \label{sec:appendix}
We give a proof of Proposition~\ref{prop:QAM=AM}. Through the quantum Yang-Baxter move \cite{KLN}, we see that it is sufficient to show the proposition for a specific reduced chain. 

\subsection{Notes on type \texorpdfstring{$C$}{C} root systems}
Assume $G = \Sp_{2n}(\BC)$. 
For $1 \le i < j \le n$, set 
\begin{equation}
(i, j) := \ve_{i} - \ve_{j}, \quad (i, \overline{j}) := \ve_{i} + \ve_{j}. 
\end{equation}
Also, for $1 \le i \le n$, set $(i, \overline{i}) := 2\ve_{i}$. 
Then we have 
\begin{equation}
\Delta^{+} = \{(i, j), (i, \overline{j}) \mid 1 \le i < j \le n \} \sqcup \{(i, \overline{i}) \mid 1\le i \le n \}. 
\end{equation}

By regarding elements $w \in W$ of the Weyl group as a signed-permutation on $[1, \overline{1}]$, we can describe a criterion of Bruhat/quantum edges in $\QBG(W)$. 
To describe the criterion, we set $\sign(k) := 1$ for $1 \le k \le n$ and $\sign(\overline{k}) := -1$ for $1 < k < n$. 
\begin{lem}[{\cite[Proposition~5.7]{Lenart}}] \label{lem:criterion}
Let $w \in W$. 
\begin{enu}
\item For $1 \le i < j \le n$, there exists a Bruhat edge $w \xrightarrow{(i, j)} ws_{(i, j)}$ in $\QBG(W)$ if and only if 
\begin{itemize}
\item $w(i) < w(j)$, and 
\item there does not exist any $i < k < j$ such that $w(i) < w(k) < w(j)$. 
\end{itemize}

\item For $1 \le i < j \le n$, there exists a quantum edge $w \xrightarrow{(i, j)} ws_{(i, j)}$ in $\QBG(W)$ if and only if 
\begin{itemize}
\item $w(i) > w(j)$, and 
\item for all $i < k < j$, we have $w(j) < w(k) < w(i)$. 
\end{itemize}

\item For $1 \le i < j \le n$, there exists an edge $w \xrightarrow{(i, \overline{j})} ws_{(i, \overline{j})}$ in $\QBG(W)$ if and only if 
\begin{itemize}
\item $w(i) < w(\overline{j})$, 
\item $\sign(w(i)) = \sign(w(\overline{j}))$, and 
\item there does not exist any $i < k < \overline{j}$ such that $w(i) < w(k) < w(\overline{j})$. 
\end{itemize}
In this case, such an edge is a Bruhat edge. 

\item For $1 \le i \le n$, there exists a Bruhat edge $w \xrightarrow{(i, \overline{i})} ws_{(i, \overline{i})}$ in $\QBG(W)$ if and only if 
\begin{itemize}
\item $w(i) < w(\overline{i})$ (or equivalently, $\sign(w(i)) = 1$), and 
\item there does not exist any $i < k < \overline{i}$ (or equivalently, $i < k \le n$) such that $w(i) < w(k) < w(\overline{i})$. 
\end{itemize}

\item For $1 \le i \le n$, there exists a quantum edge $w \xrightarrow{(i, \overline{i})} ws_{(i, \overline{i})}$ in $\QBG(W)$ if and only if 
\begin{itemize}
\item $w(i) > w(\overline{i})$ (or equivalently, $\sign(w(i)) = -1$), and 
\item for all $i < k < \overline{i}$ (or equivalently, $i < k \le n$), we have $w(\overline{i}) < w(k) < w(i)$. 
\end{itemize}
\end{enu}
\end{lem}

\subsection{Suitable chains}
Following \cite[\S 2.4]{KN}, we construct a specific chain enough to prove Proposition~\ref{prop:QAM=AM}. 

\begin{defn}[{cf. \cite[(2.1) and (2.2)]{Dyer}}]
A total order $\vtl$ on $\Delta^{+}$ is a \emph{reflection order} (or a convex order) if one of $\alpha \vtl \alpha + \beta \vtl \beta$ or $\beta \vtl \alpha + \beta \vtl \alpha$ holds for each $\alpha, \beta \in \Delta^{+}$ with $\alpha + \beta \in \Delta^{+}$. 
\end{defn}

For $\lambda \in P$, set 
\begin{align}
\Delta^{+}(\lambda)_{>0} &:= \{ \beta \in \Delta^{+} \mid \pair{\lambda}{\beta^{\vee}} > 0 \}, \\ 
\Delta^{+}(\lambda)_{=0} &:= \{ \beta \in \Delta^{+} \mid \pair{\lambda}{\beta^{\vee}} = 0 \}, \\ 
\Delta^{+}(\lambda)_{<0} &:= \{ \beta \in \Delta^{+} \mid \pair{\lambda}{\beta^{\vee}} < 0 \}. 
\end{align}
Also, for $k \in \BZ$, we set 
\begin{equation}
\Delta^{+}(\lambda)_{=k} := \{ \beta \in \Delta^{+} \mid \pair{\lambda}{\beta^{\vee}} = k \}. 
\end{equation}

For $\lambda \in P$, we denote by $\lambda_{+} \in P^{+}$ the unique element of the $W$-orbit $W\lambda$. We define $u(\lambda) \in W$ as the unique minimal element of $\{w \in W \mid w\lambda_{+} = \lambda\}$, and we denote by $w(\lambda)$ the unique maximal element of the same set as $u(\lambda)$. 
Note that if we denote by $w_{\circ}(\lambda_{+})$ the longest element of the parabolic subgroup $W_{\lambda_{+}} = \{ w \in W \mid w\lambda_{+} = \lambda_{+}\}$, then we have $w(\lambda) = u(\lambda)w_{\circ}(\lambda_{+})$. 
Then the decomposition $w_{\circ} = u(-\lambda)^{-1} u(\lambda) w_{\circ}(\lambda_{+})$ induces a reflection order $\vtl$ on $\Delta^{+}$ such that for $\alpha \in \Delta^{+}(\lambda)_{<0}$, $\beta \in \Delta^{+}(\lambda)_{=0}$, $\gamma \in \Delta^{+}(\lambda)_{>0}$, we have $\alpha \vtl \beta \vtl \gamma$ (see \cite[Remark~2.11]{KN}). 

Let $J = \{i_{1} < \cdots < i_{r}\} \subset [1, n]$ and set $\lambda = \ve_{J}$. In our setting, we have 
\begin{align}
\Delta^{+}(\lambda)_{>0} &= \{(i, j) \mid i < j, \ i \in J, \ j \notin J \} \sqcup \{(i, \overline{j}) \mid i < j, \ \text{$i \in J$ or $ j \in J$} \} \sqcup \{(i, \overline{i}) \mid i \in J\}, \\ 
\begin{split}
\Delta^{+}(\lambda)_{=0} &= \{(i, j) \mid i < j, \ i, j \in J \} \sqcup \{(i, j) \mid i < j, \ i, j \notin J \} \\ 
& \hspace{40mm} \sqcup \{(i, \overline{j}) \mid i < j, \ i, j \notin J\} \sqcup \{(i, \overline{i}) \mid i \notin J\}, 
\end{split} \\ 
\Delta^{+}(\lambda)_{<0} &= \{(i, j) \mid i < j, \ i \notin J, \ j \in J \}. 
\end{align}
To compute a reflection order $\vtl$, we use the following reduced expressions: 
\begin{align}
\begin{split}
u(-\lambda) &= (s_{i_{r}} \cdots s_{n-1} s_{n} s_{n-1} \cdots s_{1}) (s_{i_{r-1}+1} \cdots s_{n-1}s_{n}s_{n-1} \cdots s_{2}) \cdots \\ 
& \hspace{60mm} \times (s_{i_{1}+r-1} \cdots s_{n-1}s_{n}s_{n-1} \cdots s_{r}), 
\end{split} \\ 
u(\lambda) &= (s_{i_{1}-1} \cdots s_{2}s_{1})(s_{i_{2}-1} \cdots s_{3}s_{2}) \cdots (s_{i_{r}-1} \cdots s_{r+1}s_{r}), \\ 
w_{\circ}(\lambda_{+}) &= s_{1} (s_{2}s_{1}) \cdots (s_{r-1} \cdots s_{2}s_{1}) (s_{r+1} \cdots s_{n-1}s_{n})^{n-r}. 
\end{align}
Since an explicit form of $\vtl$ is so complicated, we omit the full description of $\vtl$. Instead of it, we show the explicit description of $\vtl$ on $\Delta^{+}(\lambda)_{<0}$ as follows: 
\begin{align}
& (i_{1}-1, i_{1}) \vtl \cdots \vtl (2, i_{1}) \vtl (1, i_{1}) \\ 
& \qquad \vtl (i_{2}-1, i_{2}) \vtl \cdots \vtl (i_{1}+1, i_{2}) \vtl (i_{1}-1, i_{2}) \vtl \cdots \vtl (2, i_{2}) \vtl (1, i_{2}) \\ 
& \qquad \vtl \cdots \\ 
\begin{split}
& \qquad \vtl (i_{r}-1, i_{r}) \vtl \cdots \vtl (i_{r-1}+1, i_{r}) \vtl (i_{r-1}-1, i_{r}) \vtl (i_{r-1}-2, i_{r}) \vtl \cdots \\ 
& \qquad \quad \cdots \vtl (i_{r-2}+2, i_{r}) \vtl (i_{r-2}+1, i_{r}) \vtl (i_{r-2}-1, i_{r}) \vtl (i_{r-2}-2, i_{r}) \vtl \cdots \\ 
& \qquad \quad \cdots \\ 
& \qquad \quad \cdots \vtl (i_{2}+2, i_{r}) \vtl (i_{2}+1, i_{r}) \vtl (i_{2}-1, i_{r}) \vtl (i_{2}-2, i_{r}) \vtl \cdots \\ 
& \qquad \quad \cdots \vtl (i_{1}+2, i_{r}) \vtl (i_{1}+1, i_{r}) \vtl (i_{1}-1, i_{r}) \vtl \cdots \vtl (2, i_{r}) \vtl (1, i_{r}). 
\end{split}
\end{align}
Then we can construct a reduced $\ve_{J}$-chain $\Gamma_{\vtl}(\ve_{J})$ ``suitable'' to $\vtl$ in the sense of \cite[\S 3.2]{KN} as follows: 
\begin{equation}
\Gamma_{1} = (\gamma_{1}^{(1)}, \ldots, \gamma_{p}^{(1)})
\end{equation}
with $\{\gamma_{1}^{(1)}, \ldots, \gamma_{p}^{(1)}\} = \Delta^{+}(\lambda)_{>0}$ and $\gamma_{1}^{(1)} \vtr \cdots \vtr \gamma_{p}^{(1)}$; 
\begin{equation}
\Gamma_{2} = (\gamma_{1}^{(2)}, \ldots, \gamma_{q}^{(2)})
\end{equation}
with $\{\gamma_{1}^{(2)}, \ldots, \gamma_{q}^{(2)}\} = \Delta^{+}(\lambda)_{=2}$ and $\gamma_{1}^{(2)} \vtr \cdots \vtr \gamma_{q}^{(2)}$; 
\begin{equation}
\Gamma_{3} = (-\gamma_{1}^{(3)}, \ldots, -\gamma_{r}^{(3)})
\end{equation}
with $\{\gamma_{1}^{(3)}, \ldots, \gamma_{r}^{(3)}\} = \Delta^{+}(\lambda)_{<0}$ and $\gamma_{1}^{(3)} \vtr \cdots \vtr \gamma_{r}^{(3)}$; and set 
\begin{equation}
\Gamma_{\vtl}(\ve_{J}) = (\gamma_{1}^{(1)}, \ldots, \gamma_{p}^{(1)}, \gamma_{1}^{(2)}, \ldots, \gamma_{q}^{(2)}, -\gamma_{1}^{(3)}, \ldots, -\gamma_{r}^{(3)}). 
\end{equation}

In addition, if $\Gamma = (\beta_{1}, \ldots, \beta_{r})$ is a reduced $\mu$-chain for $\mu \in P$, then we know that $\Gamma^{\ast} := (-\beta_{r}, \ldots, -\beta_{1})$ is a reduced $(-\mu)$-chain. Thus if we set $\Gamma(-\ve_{J}) := \Gamma_{\vtl}(\ve_{J})^{\ast}$, then $\Gamma(-\ve_{J})$ is a reduced $(-\ve_{J})$-chain. 
Note that $\Gamma(-\ve_{J})$ is a concatenation of three chains $\Gamma_{3}^{\ast}$, $\Gamma_{2}^{\ast}$, and $\Gamma_{1}^{\ast}$ in this order; that is., 
\begin{equation}
\Gamma(-\ve_{J}) = (\underbrace{\gamma_{r}^{(3)}, \ldots, \gamma_{1}^{(3)}}_{\Gamma_{3}^{\ast}}, \underbrace{-\gamma_{q}^{(2)}, \ldots, -\gamma_{1}^{(2)}}_{\Gamma_{2}^{\ast}}, \underbrace{-\gamma_{p}^{(1)}, \ldots, -\gamma_{1}^{(1)}}_{\Gamma_{1}^{\ast}}). 
\end{equation}

We show the following proposition, which immediately proves Proposition~\ref{prop:QAM=AM}. 
\begin{prop} \label{prop:QAM=AM_prop}
\begin{enu}
\item For $A \in \CA(e, \Gamma_{\vtl}(\ve_{J}))$, we have $\down(A) = 0$. 
\item For $A \in \CA(e, \Gamma(-\ve_{J}))$, we have $\down(A) = 0$. 
\end{enu}
\end{prop}

First, we show Proposition~\ref{prop:QAM=AM_prop}\,(1). 
For the purpose, set 
\begin{equation}
D(w) := \{ (i, j) \in [1, n] \times [1, n] \mid i < j, \ w(i) > w(j) \}. 
\end{equation}
For a directed path 
\begin{equation}
\bp: w_{0} \xrightarrow{\gamma_{1}} w_{1} \xrightarrow{\gamma_{2}} \cdots \xrightarrow{\gamma_{l}} w_{l}, 
\end{equation}
in $\QBG(W)$, set
\begin{align}
\init(\bp) &:= w_{0}, \\ 
\ed(\bp) &:= w_{l}, \\ 
\wt(\bp) &:= \sum_{\substack{1 \le j \le l \\ \text{$w_{j-1} \rightarrow w_{j}$ is a quantum edge}}} \gamma_{j}^{\vee}. 
\end{align}
Then, by using Lemma~\ref{lem:criterion}, we can show the following lemmas. 

\begin{lem} \label{lem:down=0}
Let $w \in W$ and assume that 
\begin{equation} \label{eq:descent}
D(w) \subset \{(i, j) \in [1, n] \times [1, n] \mid i < j, \ i \in J, \ j \notin J \}. 
\end{equation}

\begin{enu}
\item Let $w \xrightarrow{\gamma} w'$ be a Bruhat edge with $\gamma \in \Delta^{+}(\lambda)_{>0}$. Then we have 
\begin{equation} 
D(w') \subset \{(i, j) \in [1, n] \times [1, n] \mid i < j, \ i \in J, \ j \notin J \}. 
\end{equation}
\item Let $w = w_{0} \xrightarrow{\gamma_{1}} w_{1} \xrightarrow{\gamma_{2}} \cdots \xrightarrow{\gamma_{l}} w_{l}$ be a directed path in $\QBG(W)$ such that $\gamma_{1}, \ldots, \gamma_{l} \in \Delta^{+}(\lambda)_{<0}$ and $\gamma_{1} \vtr \cdots \vtr \gamma_{l}$. 
Then we have 
\begin{equation}
\{ \gamma_{1}, \ldots, \gamma_{l} \} \subset \{(i_{r}-1, i_{r}), \ldots, (i_{2}-1, i_{2}), (i_{1}-1, i_{1})\} \setminus \{(i_{k}, i_{l}) \mid 1 \le k < l \le r\} 
\end{equation}
and all edges in the directed path are Bruhat edges. 
\end{enu}
\end{lem}

\begin{proof}
First, we show (1). 
Note that by assumption \eqref{eq:descent}, we have $w(i) < w(j)$ for $1 \le i < j \le n$ with $i \notin J$ or $j \in J$. 
In this proof, we implicitly use Lemma~\ref{lem:criterion}. 

\paragraph{\underline{\textbf{Case~1}: $\gamma = (p, q)$ with $p < q$, $p \in J$, and $q \notin J$}} 
For $1 \le k \le n$, we see that 
\begin{equation}
w'(k) = \begin{cases} 
w(q) & \text{if $k = p$}, \\ 
w(p) & \text{if $k = q$}, \\
w(k) & \text{if $k \not= p, q$}. 
\end{cases}
\end{equation}
Let $1 \le i < j \le n$ and compare $w(i)$ and $w(j)$. 

\subparagraph{\underline{\textbf{Subcase~1.1}: $i \notin J$}}
In this case, $i \not= p$. 

If $i = q$, then $j \not= p, q$ and hence $w'(i) = w(p)$, $w'(j) = w(j)$. 
If $j \in J$, then $w(p) < w(j)$ since $p < q = i < j$ and $j \in J$. Hence $w'(i) < w'(j)$. 
Then assume that  $j \notin J$. Since $w \xrightarrow{\gamma} w'$ is a Bruhat edge, we see that $w(p) < w(q)$. Also, $w(q) < w(j)$ since $q = i < j$ and $q \notin J$. Thus $w(p) < w(j)$ and hence $w'(i) < w'(j)$. 

Assume that $i \not= q$. Then $w'(i) = w(i)$ since $i \not= p, q$. 
Let $j \in J$. If $j = p$, then $w'(j) = w(q)$. Since $i < j = p < q$ and $i \notin J$, we have $w(i) < w(q)$, and hence $w'(i) < w'(j)$. If $j \not= p$, then $w'(j) = w(j)$ since $j \not= p, q$. Then $w'(i) = w(i) < w(j) = w'(j)$ since $i < j$ and $i \notin J$. Next, let $j \notin J$. If $j = q$, then $w'(j) = w(p)$. If $i < p$, then $i \notin J$ implies $w(i) < w(p)$, and hence $w'(i) < w'(j)$. Assume that $p < i < j = q$. Since $w \xrightarrow{\gamma} w'$ is a Bruhat edge, $w(p) < w(i) < w(q)$ can not hold. On the other hand, $i < q$ and $i \notin J$ imply $w(i) < w(q)$. Hence we obtain $w(i) < w(p) < w(q)$. Therefore $w'(i) < w'(j)$. If $j \not= q$, then $w'(j) = w(j)$. Since $i \notin J$, we obtain $w'(i) = w(i) < w(j) = w'(j)$. 

\subparagraph{\underline{\textbf{Subcase~1.2}: $i \in J$ and $j \in J$}} In this case, $i, j \not= q$. 

If $i = p$, then $w'(i) = w(q)$ and $w'(j) = w(j)$. 
Assume that $i = p < j < q$. Since $w \xrightarrow{\gamma} w'$ is a Bruhat edge, $w(p) < w(j) < w(q)$ can not hold. On the other hand, $p < j$ and $j \in J$ imply $w(p) < w(j)$. Hence $w(p) < w(q) < w(j)$ and therefore $w'(i) < w'(j)$. If $q < j$, then $j \in J$ implies $w(q) < w(j)$, and hence $w'(i) < w'(j)$. 

If $j = p$, then $w'(i) = w(i)$ and $w'(j) = w(q)$. Since $i < j = p$ and $p \in J$, we see that $w(i) < w(p)$. Also, since $w \xrightarrow{\gamma} w'$ is a Bruhat edge, we have $w(p) < w(q)$. Hence $w(i) < w(q)$ and therefore $w'(i) < w'(j)$. 

If $i, j \not= p$, then $w'(i) = w(i)$ and $w'(j) = w(j)$. Since $j \in J$, we have $w'(i) = w(i) < w(j) = w'(j)$. 

\paragraph{\underline{\textbf{Case~2}: $\gamma = (p, \overline{q})$ with $1 \le p < q \le n$ and $p \in J$}}
For $1 \le k \le n$, we see that 
\begin{equation}
w'(k) = \begin{cases}
w(\overline{q}) & \text{if $j = p$}, \\ 
w(\overline{p}) & \text{if $j = q$}, \\ 
w(k) & \text{if $j \not= p, q$}. 
\end{cases}
\end{equation}

\subparagraph{\underline{\textbf{Subcase~2.1}: $j \in J$}} 

Assume that $i = p$. Then $w'(i) = w(\overline{q})$. If $j = q$, thus $w'(j) = w(\overline{p})$. Since $q = j \in J$, we have $w(p) < w(q)$. Hence $w'(i) = w(\overline{q}) < w(\overline{p}) = w'(j)$. 
If $j \not= q$, then $w'(j) = w(j)$. Since $j \in J$, we have $w(p) < w(j)$. Since $w \xrightarrow{\gamma} w'$ is a Bruhat edge, $w(p) < w(j) < w(\overline{q})$ can not hold. Hence we have $w(\overline{q}) < w(j)$ and therefore $w'(i) < w'(j)$. 

Assume that $i = q$. Then $j \not= p$. thus $w'(i) = w(\overline{p})$ and $w'(j) = w(j)$. By $j \in J$, we have $w(q) < w(j)$. Since $w \xrightarrow{\gamma} w'$ is a Bruhat edge, $w(q) < w(j) < w(\overline{p})$ (equivalently, $w(p) < w(\overline{j}) < w(\overline{q})$) can not hold. Hence $w(\overline{p}) < w(j)$ and therefore $w'(i) < w'(j)$. 

Assume that $i \not= p, q$. Then $w'(i) = w(i)$. if $j = p$, then $w'(j) = w(\overline{q})$. Since $p \in J$, we have $w(i) < w(p)$. Since $w(p) < w(\overline{q})$, we obtain $w(i) < w(\overline{q})$, and hence $w'(i) < w'(j)$. If $j = q$, then $w'(j) = w(\overline{p})$. Since $q = j \in J$, we have $w(i) < w(q)$. Since $w(p) < w(\overline{q})$, we obtain $w(i) < w(q) < w(\overline{p})$, and hence $w'(i) < w'(j)$. If $j \not= p, q$, then $w'(i) = w(i) < w(j) = w'(j)$ since $j \in J$. 

\subparagraph{\underline{\textbf{Subcase~2.2}: $i \notin J$ and $j \notin J$}} In this case, $i, j \not= p$. 

Assume that $i = q$. Then $j \not= q$ and $w'(i) = w(\overline{p})$, $w'(j) = w(j)$. Since $q = i \notin J$, we have $w(q) < w(j)$. Since $w \xrightarrow{\gamma} w'$ is a Bruhat edge, $w(q) < w(j) < w(\overline{p})$ can not hold. Hence $w(\overline{p}) < w(j)$ and therefore $w'(i) = w'(j)$. 

Assume that $j = q$. Then $w'(i) = w(i)$ and $w'(j) = w(\overline{p})$. Since $i \notin J$, we have $w(i) < w(q)$. Since $w(q) < w(\overline{p})$, we obtain $w'(i) = w(i) < w(\overline{p}) = w'(j)$. 

Assume that $i, j \not= q$. Since $i \notin J$, we have $w'(i) = w(i) < w(j) = w'(j)$. 

\paragraph{\underline{\textbf{Case~3}: $\gamma = (p, \overline{q})$ with $p < q$, $p \notin J$, and $q \in J$}}

\subparagraph{\underline{\textbf{Subcase~3.1}: $j \in J$}} In this case, $j \not= p$. 

Assume that $i = p$. Then $w'(i) = w(\overline{q})$. If $j = q$, then $w'(j) = w(\overline{p})$. Since $q \in J$, we have $w(p) < w(q)$, and hence $w'(i) = w(\overline{q}) < w(\overline{p}) = w'(j)$. 
If $j \not= q$, then $w'(j) = w(j)$. Since $j \in J$, we have $w(p) < w(j)$. Since $x \xrightarrow{\gamma} w'$ is a Bruhat edge, $w(p) < w(j) < w(\overline{q})$ can not hold. Hence $w(\overline{q}) < w(j)$ and therefore $w'(i) < w'(j)$. 

Assume that $i = q$. Then $w'(i) = w(\overline{p})$ and $w'(j) = w(j)$. Since $j \in J$, we have $w(q) < w(j)$. Since $w \xrightarrow{\gamma} w'$ is a Bruhat edge, $w(q) < w(j) < w(\overline{p})$ can not hold. Hence $w(\overline{p}) < w(j)$ and therefore $w'(i) < w'(j)$. 

Assume that $i \not= p, q$. Then $w'(i) = w(i)$. If $j = q$, then $w'(j) = w(\overline{p})$. Since $q \in J$, we have $w(i) < w(q)$. Since $w(q) < w(\overline{p})$, we obtain $w'(i) < w'(j)$. 
If $j \not= q$, then we have $w'(i) = w(i) < w(j) = w'(j)$. 

\subparagraph{\underline{\textbf{Subcase~3.2}: $i \notin J$ and $j \notin J$}} 
In this case, $i, j \not= q$. 

Assume that $i = p$. Then $w'(i) = w(\overline{q})$ and $w'(j) = w(j)$. Since $p \notin J$, we have $w(p) < w(j)$. Since $w \xrightarrow{\gamma} w'$ is a Bruhat edge, $w(p) < w(j) < w(\overline{q})$ can not hold. Hence $w(\overline{q}) < w(j)$ and therefore $w'(i) < w'(j)$. 

Assume that $j = p$. Then $w'(i) = w(i)$ and $w'(j) = w(\overline{q})$. Since $i \notin J$, we have $w(i) < w(p)$. Since $w(p) < w(\overline{q})$, we obtain $w'(i) < w'(j)$. 

Assume that $i, j \not= p$. Then we see that $w'(i) = w(i) < w(j) = w'(j)$. 

\paragraph{\underline{\textbf{Case~4}: $\gamma = (p, \overline{p})$ with $p \in J$}}
For $1 \le k \le n$, we see that 
\begin{equation}
w'(k) = \begin{cases}
w(\overline{p}) & \text{if $k = p$}, \\
w(k) & \text{if $k \not= p$}. 
\end{cases}
\end{equation}
Observe that $w'(p) > w(p)$ since $w \xrightarrow{\gamma} w'$ is a Bruhat edge. 

\subparagraph{\underline{\textbf{Subcase~4.1}: $i \notin J$}} 
In this case, we can verify clearly that $w'(i) < w'(j)$ since $w'(i) = w(i) < w(j) \le w'(j)$. 

\subparagraph{\underline{\textbf{Subcase~4.2}: $i \in J$ and $j \in J$}}

If $i \not= p$, then clearly $w'(i) < w'(j)$ by the same argument as Subcase~4.1. 
If $i = p$, then $w'(i) = w(\overline{p})$. Since $j \in J$, we have $w(p) < w(j)$. Since $w \xrightarrow{\gamma} w'$ is a Bruhat edge, $w(p) < w(j) < w(\overline{p})$ can not hold. Hence $w(\overline{p}) < w(j)$ and therefore $w'(i) < w'(j)$. 

\vspace{0.5\baselineskip}

Therefore, we conclude that if $i \notin J$ or $j \in J$, we have $w'(i) < w'(j)$. 
This proves (1). 

We hereby show (2). For this purpose, we prove the following: 
\begin{claim}
If $\gamma_{1} = (i_{t_{1}}-1, i_{t_{1}}), \ldots, \gamma_{l-1} = (i_{t_{l-1}}-1, i_{t_{l-1}})$ with $r \ge t_{1} > \cdots > t_{l-1} > 1$, then $\gamma_{l} = (i_{t_{s}}-1, i_{t_{s}})$ for some $1 \le t_{s} < t_{l-1}$ and $w_{l-1} \xrightarrow{\gamma_{l}} w_{l}$ is a Bruhat edge. 
\end{claim}

Set $w_{l-1} := w$, $w_{l} := w'$. 
We understand that $t_{0} = r+1$. 
Let $1 \le s < t_{l-1}$ and take $1 \le k < i_{s}$ wtih $k \notin J$. If $k \not= i_{s}-1$, then there exists some $k < p < i_{s}$. 
Since $k \notin J$, we have $w(k) < w(p)$. Since $i_{s} \in J$, we have $w(p) < w(i_{s})$. Hence $w(k) < w(p) < w(i_{s})$ and therefore there does not exist any Bruhat edge $w \xrightarrow{(k, i_{s})} w'$. Hence $\gamma$ must be identical to $(i_{s}-1, i_{s})$. 
In addition, since $w(i_{s}-1) < w(i_{s})$, this edge must be a Bruhat edge. This claim proves (2) of the proposition. 
\end{proof}

\begin{proof}[Proof of Proposition~\ref{prop:QAM=AM_prop}\,(1)]
Let $A \in \CA(e, \Gamma_{\vtl}(\ve_{J}))$ and for $k = 1, 2, 3$, take a subset $A^{(k)}$ of $A$ corresponding to $\Gamma_{k}$. 
Then take directed paths $\Pi_{1}, \Pi_{2}, \Pi_{3}$ in $\QBG(W)$ corresponding to $A^{(1)}$, $A^{(2)}$, $A^{(3)}$, respectively. 
Then $\Pi_{1}$ is a label-decreasing path with respect to $\vtl$ with all labels contained in $\Delta^{+}(\lambda)_{>0}$. 
Note that, $\init(\Pi_{1}) = e$. 
Since $\ed(\Pi_{1}) > e$ under the Bruhat order, we can take a shortest directed path from $e$ to $\ed(\Pi_{1})$ in $\QBG(W)$ containing only Bruhat edges. Since a label-decreasing path is shortest, we see that $\Pi_{1}$ is also shortest among all directed paths from $e$ to $\ed(\Pi_{1})$ (see \cite[Theorem~6.4]{BFP}). 
According to \cite[Lemma~1\,(2)]{Postnikov} (see also \cite[Proposition~8.1]{LNSSS}), if $\bp_{1}$ and $\bp_{2}$ are shortest directed paths such that $\init(\bp_{1}) = \init(\bp_{2})$ and $\ed(\bp_{1}) = \ed(\bp_{2})$, then $\wt(\bp_{1}) = \wt(\bp_{2})$. 
Hence we have $\wt(\Pi_{1}) = 0$. 
Also, by Lemma~\ref{lem:criterion}\,(3), we have $\wt(\Pi_{2}) = 0$. 
By Lemma~\ref{lem:down=0}\,(1), $w := \ed(\Pi_{2})$ satisfies \eqref{eq:descent}. Then by Lemma~\ref{lem:down=0}\,(2), we have $\wt(\Pi_{3}) = 0$. 
Hence $\down(A) = \wt(\Pi_{1}) + \wt(\Pi_{2}) + \wt(\Pi_{3}) = 0$, as desired. 
\end{proof}

To show Proposition~\ref{prop:QAM=AM_prop}\,(2), we see the following lemma, which can be easily verified from the following fact:  
\begin{equation}
\Delta^{+}(\lambda)_{<0} \subset \{(i, j) \mid 1 \le i < j \le n\}. 
\end{equation} 
\begin{lem}
Let $e \xrightarrow{\gamma_{1}} \cdots \xrightarrow{\gamma_{l}} w$ be a directed path in $\QBG(W)$ with $\gamma_{1}, \ldots, \gamma_{l} \in \Delta(\lambda)_{<0}$. Then we have $w(k) \in \{ 1, \ldots, n \}$ for $1 \le k \le n$. 
\end{lem}

\begin{cor} \label{cor:no_edge}
Let $e \xrightarrow{\gamma_{1}} \cdots \xrightarrow{\gamma_{l}} w$ be a directed path in $\QBG(W)$ such that 
\begin{itemize}
\item all edges are Bruhat edges, and 
\item $\gamma_{1}, \ldots, \gamma_{l} \in \Delta^{+}(\lambda)_{<0}$. 
\end{itemize}
Then there does not exist any edge $w \xrightarrow{\gamma} w'$ in $\QBG(W)$ such that $\gamma \in \Delta^{+}(\lambda)_{=2}$. 
\end{cor}
\begin{proof}
According to Lemma~\ref{lem:criterion}\,(3), this follows from the fact $\sign(w(i)) = 1 \not= -1 = \sign(w(\overline{j}))$ for all $1 \le i < j \le n$. 
\end{proof}

\begin{proof}[Proof of Proposition~\ref{prop:QAM=AM_prop}\,(2)]
Let $A \in \CA(e, \Gamma(-\ve_{J}))$ and for $k = 1, 2, 3$, take a subset $A^{(k)}$ of $A$ corresponding to $\Gamma_{k}^{\ast}$. 
For $k = 1, 2, 3$, let $\Pi_{k}^{\ast}$ be a directed path in $\QBG(W)$ corresponding to $A^{(k)}$. 
Since $\ed(\Pi_{3}) > e$ under the Bruhat order and $\Pi_{3}$ is a label-increasing path in $\QBG(W)$ with respect to $\vtl$, we see that $\wt(\Pi_{3}^{\ast}) = 0$ by the same reason for the identity $\wt(\Pi_{1}) = 0$ in the proof of Proposition~\ref{prop:QAM=AM_prop}\,(1). Then by Corollary~\ref{cor:no_edge}, we see that $\Pi_{2}^{\ast}$ is a trivial path (of length $0$). 
Therefore, the directed path $\Pi(e, A)$ corresponding to $A$ is a concatenation of $\Pi_{3}^{\ast}$ and $\Pi_{1}^{\ast}$ in this order, and this is a label-increasing path in $\QBG(W)$ with respect to $\vtl$. Since $\ed(A) > e$ under the Bruhat order, we conclude that $\down(A) = \wt(\Pi(e, A)) = 0$, as desired. 
\end{proof}

\end{document}